\documentclass[a4paper,twoside,11pt]{article}
\usepackage[hiresbb]{graphicx}
\usepackage{amsmath, amssymb,amsthm, amsfonts}
\usepackage{color, subcaption, enumitem, hyperref}

\usepackage{newtxtext}
\usepackage[varvw]{newtxmath}

\usepackage{bm}
\usepackage[mathscr]{eucal}
\usepackage{authblk}
\usepackage{ascmac}
\usepackage{wrapfig}
\usepackage[english]{babel}
\usepackage{scirep-A4} 
\usepackage[]{multicol}
\usepackage{here}
\usepackage{empheq}
\usepackage{listings}

 \DeclareSymbolFont{matha}{OML}{txmi}{m}{it}

\DeclareMathOperator*{\essinf}{ess\,inf}

\newtheorem{theorem}{Theorem}[section]

\newtheorem{remark}[theorem]{Remark}
\theoremstyle{definition}

\newcommand{\DIV}{\operatorname{div}}

\def\N{{\mathbb N}}
\def\R{{\mathbb R}}

\def\DS{\displaystyle}

\makeatletter
 \def\fgcaption{\def\@captype{figure}\caption}
 \makeatother

\title{\textbf{
Energy-consistent dynamic fracture phase field models: unilateral constraints and finite element simulations
}}
\author{
Md Mamun {\sc Miah}$^{1,2}$\thanks
{{\it Corresponding author} \quad 
Email: \ mamun0954@gmail.com}~, 
Ryuhei {\sc Wakida}$^{1}$, \ and \ Masato {\sc Kimura}$^{3}$
\vspace{1.5mm}
\\
{\small 
$^1$ Division of Mathematical and Physical Sciences, Kanazawa University,Kakuma,
}\\
{\small 
Kanazawa, 920-1192, Japan
}\\
{\small 
$^2$ Department of Mathematics, Khulna University of Engineering and Technology,
}\\
{\small 
Khulna-9203, Khulna, Bangladesh
}\\
{\small 
$^3$ Faculty of Mathematics and Physics, Kanazawa University, Kakuma,
}\\
{\small 
Kanazawa, 920-1192, Japan
}\\
}

\date{\small{
}} 

\begin{document}
\maketitle
\makeatletter

\begin{minipage}{.94\linewidth}
{\bf Abstract}\quad
Phase field models have emerged as a powerful and flexible framework for simulating complex interface-driven phenomena across a wide range of scientific and engineering applications. In fracture mechanics, the phase field approach—formulated as a gradient flow of the Griffith fracture energy with Ambrosio-Tortorelli regularization—has gained significant attention for its ability to capture complex crack topologies.
In this study, we propose a dynamic fracture phase field model (DF-PFM) based on the elastodynamic wave equation. We further extend this framework by incorporating a unilateral contact condition, yielding a refined model suitable for simulating fault rupture under high pressure. For both models, we formally derive energy dissipation identities under mixed boundary conditions, 
providing insights into the energetic structure
of the formulations.
To validate the proposed approach, we conduct numerical experiments using linear implicit time discretization and finite element methods. Our simulations demonstrate that the unilateral contact condition is essential for accurately capturing shear-dominated crack propagation and preventing non-physical interpenetration, especially under high-compression loading scenarios relevant to seismic faulting.
\medskip

%
\end{minipage}
\medskip

\normalsize
\thispagestyle{firstheadings}
\pagestyle{myheadings}

\markboth{{\footnotesize\rm
Energy-consistent dynamic fracture models with contact
}}
{{\footnotesize\rm 
Md Mamun {\sc Miah}, Ryuhei {\sc Wakida} and Masato {\sc Kimura}
}}
\bigskip

\section{Introduction}
Fracture phenomena under extreme loading conditions are central to many engineering applications, including hydraulic fracturing in petroleum engineering, crash safety in automotive and aerospace industries, structural integrity assessment, nondestructive testing, health monitoring, and seismic fault modeling. Among these, the simulation and prediction of dynamic crack propagation remains a fundamental challenge in solid mechanics, despite decades of intense research. In particular, understanding and reproducing fault-type earthquake ruptures—those that occur within tectonic plates—has become an urgent topic in modern geoscience and engineering.

Numerical simulation techniques are indispensable tools for analyzing crack propagation due to their versatility and cost-effectiveness. In this context, the phase field model (PFM) has emerged as a powerful approach for fracture mechanics. The fracture phase field model (F-PFM), originally introduced as a regularization of Griffith's fracture energy using the Ambrosio–Tortorelli functional, allows for a smooth approximation of cracks without explicit tracking of discontinuities. This model handles complex crack topologies naturally, simplifies numerical implementation, and enables autonomous crack path prediction.

The theoretical underpinnings of phase field models can be traced back to the Cahn–Hilliard theory, which governs phase separation in binary coexisting phases~\cite{Cahn_1958}. 
Since then, the methodology has been adapted for a wide range of interfacial phenomena~\cite{Allen_1972, Caginalp_1989, Kobayashi_1993}. In the field of fracture mechanics, Bourdin et al.~\cite{Bourdin_2000} and Karma et al.~\cite{Karma_2001} were among the first to apply phase field methods to variational fracture modeling. Later, Takaishi and Kimura~\cite{TK_2009, Kimura2_2021} derived an irreversible F-PFM using a gradient flow of the Francfort–Marigo energy~\cite{Francfort_1998}, incorporating irreversibility via a unidirectional gradient flow~\cite{Akagi_2019}.

Since then, the method has been extended to various complex scenarios such as fatigue~\cite{Alessi_2018}, 
anisotropic~\cite{Teichtmeister_2017}, 
three-dimensional micro~\cite{Nguyen_2016}, 
concrete damage~\cite{Comi_2001},
thermal~\cite{Alfat_2024, Mandal_2021, Miehe_2015, Zhang_2023}, 
(visco)plastic~\cite{Arriaga_2017,Guo_2005},
poroelastic~\cite{Miehe_2016},
hyperelastic~\cite{Xing_2023},
and dynamic~\cite{Borden_2012, Schluter_2014} 
fracture problems.
Recent attempts to further extend the fracture phase field model and their analysis include the following: 
\cite{Freddi_2010, Henry_2004, Leng_2024, Li_2023, Liu_2023a, Rudshaug_2024, Si_2024}

Despite its many advantages, standard F-PFM formulations can exhibit unphysical behavior in compression-dominated scenarios due to their inability to account for contact between crack faces. This issue is particularly relevant in dynamic rupture problems such as those found in fault-type earthquakes. To address this, unilateral contact conditions have been incorporated into phase field models~\cite{Amor_2009,Liu_2021, Nguyen_2019, Nguyen_2020}, which enforce non-penetration constraints and lead to more realistic simulations under mixed-mode and compressive loading conditions.

In the present study, we investigate dynamic fracture propagation in the context of fault-type earthquakes using a dynamic extension of F-PFM (DF-PFM) based on the elastodynamic wave equation. We further propose and implement a new model that incorporates a unilateral contact condition into the DF-PFM framework. For both models, we derive energy dissipation identities under appropriate boundary conditions to ensure theoretical consistency. 

We also propose linear implicit time-discrete schemes for both formulations, derive their corresponding weak forms, and perform numerical experiments using finite element methods. Our simulations focus on a high-pressure fault zone with an inclined initial crack subjected to an incident P-wave. We demonstrate that, while the original DF-PFM exhibits unrealistic kink-type crack patterns due to negative opening displacements, the proposed DF-PFM with a unilateral contact condition produces physically plausible shear-dominated rupture propagation.

To our knowledge, this is the first study to systematically address the limitations of irreversible-type dynamic fracture PFM in dynamic compression scenarios and to offer a validated remedy via contact-aware modeling.

The structure of this paper is as follows: Section~2 introduces the energy-consistent dynamic fracture phase field model. Section~3 extends the model to include a unilateral contact condition. Section~4 formulates the time-discrete schemes and corresponding weak forms. Section~5 presents two-dimensional finite element simulations of fracture dynamics. Section~6 discusses the key findings, and Section~7 concludes the paper.

\section{The dynamic fracture phase field model}
In this section, we introduce the dynamic fracture phase field model, which extends the fracture phase field model proposed in \cite{Kimura2_2021, TK_2009} by adding an inertia term to the equation. We then examine the nature of its energy balance.

Let $d=2$ or $3$. Consider a  bounded Lipschitz domain $\Omega \subset \mathbb{R}^n$ with boundary $\partial\Omega = \Gamma = \Gamma_D \cup \Gamma_N$, where $\Gamma_D$ and $\Gamma_N$ are disjoint measurable subsets of $\Gamma$.
We assume that the interface $\overline{\Gamma_D} \cap \overline{\Gamma_N}$ is sufficiently regular. The outward normal vector $\nu \in \R^d$ is defined on $\Gamma$. 
We impose Dirichlet boundary conditions on $\Gamma_D$ and Neumann boundary conditions on $\Gamma_N$ for the displacement.

For two tensors (or matrices) $\xi=(\xi_{ij})$ and $\eta=(\eta_{ij})\in \R^{d\times d}$,
We define the inner product as
$\xi:\eta\coloneqq \xi_{ij}\eta_{ij}=\text{tr}(\xi^T\eta)$, where we use the Einstein summation convention. 
The norm of $\xi\in \R^{d\times d}$
is defined by 
$|\xi| \coloneqq \sqrt{\xi :\xi }$. 

The domain $\Omega$ models a brittle elastic solid.
Let $x \in \overline\Omega$ 
denote the position and $t\in [0,T]$ denote the time.
For a displacement field $u=(u_i)$~:~$\overline\Omega \to \R^d$, 
the strain and stress tensors are defined as follows
\begin{align}
   & e[u] = (e_{ij}[u]) \in \mathbb{R}^{d\times d}_{\text{sym}}, \quad \text{where}~~~ e_{ij}[u] \coloneqq \frac{1}{2}\left(\frac{\partial u_i}{\partial x_j} + \frac{\partial u_j}{\partial x_i} \right), \notag\\
   & \sigma[u] =(\sigma_{ij}[u])\coloneqq  C e[u]\in\mathbb{R}^{d\times d}_{\text{sym}}, \quad \text{i.e.,} \quad \sigma_{ij}[u] \coloneqq c_{ijkl} e_{kl}[u]. \notag
\end{align}

The elasticity tensor $C(x) = (c_{ijkl}(x))$ satisfies the symmetry conditions:
\begin{align*}
c_{ijkl}(x)=c_{klij}(x)=c_{jikl}(x)\quad (x\in\Omega)
\end{align*}
and the coercivity condition
\begin{align}\label{coercive}
~^\exists c_*>0\quad \text{s.t.}~c_{ijkl}(x)\xi_{ij} \xi_{kl} \ge c_* |\xi|^2\quad \left(x\in\Omega,~\xi=(\xi_{ij})\in \mathbb{R}^{d\times d}_{\text{sym}}\right).
\end{align}

In particular, for homogeneous isotropic materials, the stress tensor can be expressed using the Lamé parameters $\lambda$ and $\mu$ as
\begin{align}
\sigma[u] = \lambda (\operatorname{div} u) I + 2\mu e[u],\label{stress}
\end{align}
where $I\in \R^{d\times d}_\text{sym}$ is the identity matrix.
The coercivity condition \eqref{coercive} holds with
$c_*= 2\mu + d \min (\lambda, 0)$ provided $\mu >0$ and $d\lambda +2\mu >0$. 

The elastic energy density for the displacement $u$ is defined by 
\begin{align}\label{Wu}
W(u) \coloneqq \sigma[u] : e[u]\quad (u\in H^1(\Omega;\R^d)),
\end{align} 
and naturally appears in our fracture model as the driving force term for failure due to stress concentration at crack tips.

In this paper, a volume force $f(x,t)\in \R^d$ in $\Omega$, a surface force $q(x,t)\in \R^d$ on $\Gamma_N$, and a prescribed displacement $g(x,t)\in\R^d$ on $\Gamma_D$, will be called the external forces.
Material parameters are defined as follows: $\rho > 0$ is the material density; $\gamma_*=\gamma_*(x)$ is the fracture energy (the critical energy release rate $G_\text{c}$ according to Griffith's criterion) satisfying
\[
\gamma_*\in L^\infty(\Omega),\quad
\essinf_\Omega \gamma_* >0,
\]
and $\epsilon > 0$ characterizes the small internal length scale. 

We consider the following initial-boundary value problem coupling elastic wave propagation and crack evolution:
\begin{subequations}
\begin{empheq}
[left=\empheqlbrace]{align}
& \hspace{5pt} 
\rho\frac{\partial^2 u}{{\partial t}^2} = \DIV \sigma [u,z] + f(x,t)  &&  \text{in}~ {\Omega}\times[0,{T}],  \label{DF1}\\
& \hspace{5pt} 
\alpha \frac{\partial z}{\partial {t}} = \left({\epsilon}~\operatorname{div}\left({\gamma_{*}} {\nabla} z\right) - \frac{{\gamma_{*}}}{{\epsilon}}z + (1-z) W(u)\right)_{+} &&\text{in}~ {\Omega}\times[0,{T}], \label{DF2}\\
& \hspace{5pt} 
\sigma[u,z]\nu = q(x,t) && \text{on}~ {\Gamma_N}\times[0,{T}],\label{DF3}\\
& \hspace{5pt} 
u = g(x,t) && \text{on}~ {\Gamma_D}\times[0,{T}],\label{DF4}\\
& \hspace{5pt} 
\frac{\partial z}{\partial \nu} = 0 && \text{on}~ {\Gamma}\times[0,{T}],\label{DF5}\\
& \hspace{5pt} 
\left.\left(u, \frac{\partial u}{\partial t}\right)\right|_{t=0} = (u_0(x), v_0(x)) && \text{in}~ {\Omega},\label{DF6}\\
& \hspace{5pt}
z|_{t=0} = z_0(x) && \text{in}~ {\Omega}.\label{DF7}
\end{empheq}\label{DF}
\end{subequations}

We call \eqref{DF} dynamic fracture phase field model (DF-PFM) in this paper.
The solution of the DF-PFM is given by the displacement $u(x,t)\in\R^d$ and damage variable $z(x,t)\in [0,1]$, also called the phase field variables. The damage variable $z$ represents the crack phase: $z \approx 0$ indicates intact material, $z \in (0,1)$ partial damage, and $z \approx 1$ full fracture. Initial data are given by $u_0(x)$ and $v_0(x)$ for $u$, and $z_0(x)$ for $z$.

Following \cite{Kimura2_2021}, an important aspect of the irreversible fracture phase field model (F-PFM) is the damaged elastic modulus defined as $\tilde{C}\coloneqq (1-z)^2C$. The associate stress tensor is
\begin{align}\label{suz}
\sigma[u,z] \coloneqq  (1-z)^2 \sigma[u]\quad 
(u\in H^1(\Omega;\R^d),~z\in L^\infty(\Omega)).
\end{align}
Equation \eqref{DF1} thus describes an elastic wave equation with damaged stress tensor $\sigma[u,z]$, balancing inertia, internal tension, and external forces.

Since cracks are irreparable, the positive part operator $(\cdot)_+= \max(\cdot, 0)$ is applied on the right-hand side of \eqref{DF2}.
This is another key issue of the irreversible F-PFM, that is theoretically supported by the mathematical studies of the irreversible gradient flow (unidirectional diffusion equation) in \cite{Akagi_2019,K-N2021}, and also by the physical implication of the time constant $\alpha >0$
 \cite{KTT_2024}.

If we set $\rho = 0$ and neglect the initial condition for $u$ \eqref{DF6} in the dynamic F-PFM \eqref{DF}, we recover the original F-PFM proposed in \cite{Kimura2_2021, TK_2009}. In that model, the total energy is defined as the sum of the elastic and interfacial energies, regularized by the Ambrosio-Tortorelli approach \cite{Ambrosio_1992}, and the F-PFM is derived as an irreversible gradient flow of this regularized energy. As shown in \cite{Kimura2_2021, TK_2009}, crack irreparability and energy dissipation equality hold simultaneously.

For simplicity, we denote the partial derivative with respect to time by $\dot{u} \coloneqq  \frac{\partial u}{\partial t}$, throughout this paper, except for $\dot{F}$ below, which has a different meaning.

For the DF-PFM \eqref{DF}, the total energy is defined by
\begin{align*}
E(t, u, v, z) &\coloneqq  E_{\rm el}(t, u, z) + E_{\rm ki}(v) + E_{\rm s}(z)\\
&(t\in [0,T],~u\in H^1(\Omega;\R^d),
~v\in L^2(\Omega;\R^d),~z\in H^1(\Omega)\cap L^\infty(\Omega)),
\end{align*}
where 
\begin{align}
E_{\rm el}(t, u, z) &\coloneqq \frac{1}{2} \int_\Omega \sigma[u,z] : e[u]\, dx - \int_\Omega f(t) \cdot u\, dx - \int_{\Gamma_N} q(t) \cdot u\, ds \label{oee}\\
E_{\rm ki}(v) &\coloneqq \frac{\rho}{2} \int_\Omega |v|^2\, dx \label{oke}\\
E_{\rm s}(z) &\coloneqq \frac{1}{2} \int_\Omega \gamma_* \left( \epsilon |\nabla z|^2 + \frac{z^2}{\epsilon} \right) dx \label{ose}
\end{align}
represent the elastic energy including external forces, kinetic energy, and surface energy, respectively.

Following \cite{Kimura2_2021, KTT_2024}, the energy injection rate due to external forces is defined by
\begin{equation}
    \dot F(t, u, z) \coloneqq \int_{\Gamma_D} \dot g(t) \cdot (\sigma[u,z] \nu)\, ds - \int_\Omega \dot f(t) \cdot u\, dx - \int_{\Gamma_N} \dot q(t) \cdot u\, ds \label{oei}.
\end{equation}

We emphasize that the following identity is derived under the assumption that all functions involved are sufficiently smooth to allow formal manipulations, such as differentiation and integration by parts. While a rigorous justification (e.g., in the weak solution framework) would require significant technical development, we restrict ourselves to the formal level here, as it suffices to illustrate the key dissipation structure of the model.
\begin{theorem}[Formal energy dissipation identity I]\label{Theorem1:a}
Assume that $(u,z)$ is a sufficiently smooth solution to \eqref{DF} so that all terms below are well-defined. Then the following (formal) energy dissipation identity holds:
\begin{align}\label{Eq:Theorem2:a}
\frac{d}{dt}E(t,u(t),\dot{u}(t),z(t)) = -\alpha\int_{\Omega} |\dot z(t)|^{2}\, dx+\dot F(t,u(t),z(t)) \quad \text{a.e. } t\in (0,T).
\end{align}
\end{theorem}
\proof 
First, we note that, in general, “a.e $t\in (0,T)$” in equation \eqref{Eq:Theorem2:a} can not be written as $t\in [0,T]$. This is due to the non-smooth nonlinearity in \eqref{DF2}. We leave a more detailed discussion to \cite{Akagi_2019,OMAM_2025} and omit a.e. for notational simplicity in the following.

We also remark on the following identity:
\begin{align*}
\frac{\partial}{\partial t} 
\left(\frac{1}{2} \sigma[u(t),z(t)] : e[u(t)] \right)
= \sigma[u,z] : e[\dot u] - (1-z) \dot z \, W(u) , 
\end{align*}
which is derived from the symmetry of the elasticity modulus $C(x)$.
This implies
\begin{align}
\frac{d}{dt}\left(\frac{1}{2}\int_\Omega \sigma [u,z]:e[u] \, dx\right)=\int_\Omega (\sigma [u,z]:e[\dot u]-(1-z)\dot z W(u))\, dx. 
\label{q1}
\end{align}

Differentiating \eqref{oee} and using  \eqref{q1},
\begin{align}
    \frac{d}{dt} E_{\rm el}(t,u(t), z(t))=&\int_\Omega \sigma[u,z]:e[\dot{u}]\, dx-\int_\Omega (1-z)W(u)\dot{z}\, dx\notag\\ &-\int_\Omega (\dot{f}\cdot u+f\cdot \dot{u})\, dx -\int_{\Gamma_N} (\dot{q}\cdot u+q\cdot \dot{u})\, ds. \label{dee}
\end{align} 
Now, using the integration by parts,
\begin{align}
\int_\Omega &\sigma[u,z]:e[\dot{u}]\, dx\notag\\
&=  
\int_\Omega \sigma[u,z]:\nabla^T \dot{u}\, dx\notag \\
&=
\int_\Gamma (\sigma[u,z]\nu)\cdot \dot{u}\, ds-\int_\Omega \DIV(\sigma[u,z])\cdot \dot{u}\, dx \notag \\
&=
\int_{\Gamma_D}(\sigma[u,z]\nu)\cdot \dot{g}\ ds+\int_{\Gamma_N}q\cdot \dot{u}\, ds-\int_\Omega (\rho\ddot{u}-f)\cdot \dot{u}\, dx.~~ (\mbox{using} \ \eqref{DF})  \label{ibp}
\end{align} 

Differentiating \eqref{oke} and \eqref{ose}, we obtain
\begin{align}\label{dke}
    \frac{d}{dt}E_{\rm ki}(\dot{u}(t))&=\frac{d}{dt} \left (\frac{\rho}{2}\int_\Omega |\dot u|^2\, dx\right)=\rho \int_\Omega \ddot {u} \cdot \dot{u}\, dx, 
\end{align} 
and 
\begin{align}
    \frac{d}{dt}E_{\rm s}(z(t))&= \frac{d}{dt}\left(\frac{1}{2}\int_\Omega \gamma_* (\epsilon|\nabla z|^2+\frac{z^2}{\epsilon})\, dx\right) \notag\\
    &=\int_\Omega \gamma_* \left(\epsilon \nabla z \cdot \nabla \dot{z}+\frac{1}{\epsilon} z\ \dot{z}\right)\, dx \notag \\
    & 
    =\int_\Gamma \gamma_* \epsilon \frac{\partial z}{\partial \nu}\dot{z}\,dS
    -\int_\Omega \left(\epsilon \DIV(\gamma_*\nabla z)-\gamma_*\frac{z}{\epsilon} \right) \dot{z}\, dx \quad (\text{int. by parts})
    \notag\\
    &
    = -\int_\Omega \left(\epsilon \DIV(\gamma_*\nabla z)-\gamma_*\frac{z}{\epsilon} \right) \dot{z}\, dx. \quad (\text{using }\eqref{DF5})\label{dse}
\end{align}
Taking the sum of \eqref{dee}, \eqref{dke}, and \eqref{dse}, and using \eqref{oei} and \eqref{ibp},
we obtain 
\begin{align*}
\frac{d}{dt}E(t,u(t),\dot{u}(t),z(t)) 
&=
 -\int_\Omega \left(\epsilon \DIV(\gamma_*\nabla z)-\gamma_*\frac{z}{\epsilon} +(1-z)W(u)\right) \dot{z}\, dx+\dot F(t,u(t),z(t)) \\
&=-\alpha\int_{\Omega} |\dot z(t)|^{2}\, dx
+\dot F(t,u(t),z(t)),
\end{align*}
where we used the following implication:
\begin{align}\label{alpha}
\alpha >0,~a,b\in\R,~\alpha a =(b)_+
\Longrightarrow
ab =\alpha a^2.
\end{align}
Hence, we conclude
\eqref{Eq:Theorem2:a}. 
\qed

\section{The dynamic fracture PFM with a unilateral contact condition}

In the context of simulating fracture phenomena, such as subsurface fault rupture, where shear stress is the primary factor driving crack propagation under high pressure, it is imperative to note that the absence of the unilateral contact condition can result in an unrealistic fracture, characterized by negative aperture displacement. Consequently, a methodology for imposing unilateral contact conditions on phase field models for fracture has been proposed by Amor et al. \cite{Amor_2009}.

The objectives of this section are twofold: first, to incorporate the unilateral condition into the dynamic F-PFM \eqref{DF}, and second, to demonstrate the validity of the energy dissipation identity even with incorporating the unilateral condition.

We suppose the homogeneous isotropic elasticity as \eqref{stress}.
We introduce the deviatoric part of the strain tensor $e_D[u]$ and bulk modulus $\lambda^*$ as 
\begin{align*}
e_D[u]\coloneqq e[u] - \frac{1}{d}(\DIV u)I,\quad
\lambda^*\coloneqq \lambda +\frac{2\mu}{d}.
\end{align*}
Then, the stress tensor is decomposed into the spherical and deviatoric parts:
\begin{align}\notag 
\sigma[u] = \lambda^*(\DIV u) I + 2\mu e_D[u].
\end{align}
Applying the Jordan decomposition to $\DIV u =(\DIV u)_+-(\DIV u)_-$, we obtain
\begin{align}\label{sigma3}
\sigma[u] = \lambda^*(\DIV u)_+ I  + 2\mu e_D[u]-\lambda^*(\DIV u)_- I=
\sigma_+[u] -\sigma_-[u],
\end{align}
where $(a)_+ = \max(a, 0)$ and $(a)_- = \max(-a, 0)$, and 
\begin{align*}
\sigma_+[u] \coloneqq \lambda^*(\DIV u)_+ I + 2\mu e_D[u],
\quad
\sigma_-[u] \coloneqq \lambda^*(\DIV u)_- I.
\end{align*}
The stress tensors $\sigma_\pm [u]$ represent the stresses related to the expansion and shear,
and the compression, respectively.

Then, unlike the DF-PFM \eqref{DF}, the elastic energy density $W(u)$ defined by \eqref{Wu} that drives crack propagation is replaced by
\begin{align*}
W_{+}(u)\coloneqq\sigma_{+}[u]:e[u]\quad (u\in H^1(\Omega;\R^d)),
\end{align*}
that is the elastic energy density with respect to the expansion and shear deformation. On the other hand, the stress tensor $\sigma [u,z]$ defined by \eqref{suz} is replaced by the following term:
\begin{align*}
\sigma^\dagger [u,z]\coloneqq (1-z)^2\sigma_+[u]-\sigma_- [u]\quad 
(u\in H^1(\Omega;\R^d),~z\in L^\infty(\Omega)).
\end{align*}
which considers damage only in the expansion and shear part. 

The DF-PFM with a unilateral contract condition is given by DF-PFM \eqref{DF} with \eqref{DF1}-\eqref{DF3} replaced by the following equations:  
\begin{subequations}
\begin{empheq}
[left=\empheqlbrace]{align}
& \hspace{5pt} \rho\frac{\partial^2 u}{{\partial t}^2}=\DIV \sigma^{\dagger}[{u,z}]+f(x,t)  &&  \mbox{in}~ {\Omega}\times[0,{T}],  \label{DFU1}\\
& \hspace{5pt} 
{\alpha} \frac{\partial z}{\partial {t}} = \left({\epsilon}~\DIV \left({\gamma_{*}} {\nabla} z\right) - \frac{{\gamma_{*}}}{{\epsilon}}z + (1-z) W_+(u)\right)_{+} && \mbox{in}~ {\Omega}\times[0,{T}], \label{DFU2} \\
& \hspace{5pt} \sigma^{\dagger}[u,z]\nu=q(x,t) && \mbox{on}~ {\Gamma_N}\times[0,{T}].\label{DFU3}
\end{empheq}\label{DFU}
\end{subequations}

For the DF-PFM with a unilateral contact condition \eqref{DFU} and \eqref{DF4}-\eqref{DF7}, the total energy is defined as 
\begin{align*}
E^\dagger(t, u, v, z) &\coloneqq  E_{\rm el}^\dagger (t, u, z) + E_{\rm ki}(v) + E_{\rm s}(z)\\
&(t\in [0,T],~u\in H^1(\Omega;\R^d),
~v\in L^2(\Omega;\R^d),~z\in H^1(\Omega)\cap L^\infty(\Omega)),
\end{align*}
where 
\begin{align}\label{eeld}
E_{\rm el}^\dagger (t, u, z) &\coloneqq \frac{1}{2} \int_\Omega \sigma^\dagger[u,z] : e[u]\, dx - \int_\Omega f(t) \cdot u\, dx - \int_{\Gamma_N} q(t) \cdot u\, ds.
\end{align}

The energy injection rate due to external forces is similarly defined as:
\begin{equation}\label{dFd}
    \dot F^\dagger (t, u, z) \coloneqq \int_{\Gamma_D} \dot g(t) \cdot (\sigma^\dagger [u,z] \nu)\, ds - \int_\Omega \dot f(t) \cdot u\, dx - \int_{\Gamma_N} \dot q(t) \cdot u\, ds.
\end{equation}

As with Theorem~\ref{Theorem1:a}, even with a unilateral contact condition, the following identity is formally derived under the assumption of sufficient regularity of the solution.
\begin{remark}
{\rm The following result is derived formally. Due to the presence of non-smooth nonlinear terms such as $(\DIV u)_\pm$ in the model, we assume only that the solution $(u,z)$ is regular enough to justify each step of the derivation (such as time differentiation and integration by parts). A rigorous mathematical justification would require a detailed functional setting, including a weak formulation and regularity theory, which is beyond the scope of this paper.}
\end{remark}
\begin{theorem}[Formal energy dissipation identity II]\label{Theorem1:a2}
Assume that $(u,z)$ is a sufficiently smooth solution to \eqref{DFU} and \eqref{DF4}-\eqref{DF7} so that all terms below are well-defined. Then the following (formal) energy dissipation identity holds:
\begin{eqnarray}
\frac{d}{dt}E^\dagger(u(t),z(t),t) = -\alpha\int_{\Omega} |\dot z(t)|^{2}\, dx+\dot F^\dagger(u(t),z(t),t) ~~\mbox{a.e.} ~t\in (0,T). \label{Eq:Theorem2:a2}
\end{eqnarray}
\end{theorem}
\proof 
The proof is analogous to that of Theorem~\ref{Theorem1:a}. It is important to note that, in accordance with (2.11), the subsequent equality
\begin{align}
    \frac{d}{dt}\left(\frac{1}{2}\int_\Omega \sigma^\dagger[u,z]:e[u] \, dx\right)=\int_\Omega (\sigma^\dagger[u,z]:e[\dot u]-(1-z)\dot z W_+(u))\, dx \label{q2}
\end{align}
 is valid for a.e. $t\in (0,T)$. 
This conclusion has been demonstrated in \cite{OMAM_2025} under more stringent assumptions regarding the regularity of the solution.

Differentiating \eqref{eeld} and using  \eqref{q2},
\begin{align}
\frac{d}{dt} E_{\rm el}^\dagger (t,u(t), z(t))=&\int_\Omega \sigma^\dagger [u,z]:e[\dot{u}]\, dx-\int_\Omega (1-z)W_+(u)\dot{z}\, dx\notag\\ &-\int_\Omega (\dot{f}\cdot u+f\cdot \dot{u})\, dx -\int_{\Gamma_N} (\dot{q}\cdot u+q\cdot \dot{u})\, ds. \label{dee2}
\end{align} 
Similar to the calculation in \eqref{ibp}, 
we obtain
\begin{align}
\int_\Omega &\sigma^\dagger [u,z]:e[\dot{u}]\, dx\notag\\
&=
\int_{\Gamma_D}(\sigma^\dagger [u,z]\nu)\cdot \dot{g}\ ds+\int_{\Gamma_N}q\cdot \dot{u}\, ds-\int_\Omega (\rho\ddot{u}-f)\cdot \dot{u}\, dx.~~ (\mbox{using} \ \eqref{DFU})  \label{ibp2}
\end{align} 

Hence, \eqref{Eq:Theorem2:a2} also follows 
from \eqref{dke}, \eqref{dse}, \eqref{alpha},
\eqref{dFd},  \eqref{dee2}, and \eqref{ibp2}.
\qed

\section{Time discretization}
In this section, we discuss the time discretization of the dynamic fracture models and their weak formulations in the subsequent two subsections.

For $m\in\N$, $m\ge 2$, we define $\tau\coloneqq T/m$.
For $k=0,1,2,3,...m$, setting $t_k\coloneqq k \tau$,
we consider a time discrete approximation 
$u^k\approx u(\cdot,t_k)$ and 
$z^k\approx z(\cdot,t_k)$.
For the initial conditions, we set
\begin{align}\label{uz0}
u^0\coloneqq u_0,\quad
u^1\coloneqq u_0+\tau v_0,\quad
z^0\coloneqq z_0.
\end{align}

For simplicity of description, this section assumes $f=q=0$, and $g\in C^0([0,T];H^1(\Omega;\R^d))$.
For $b\in H^1(\Omega;\R^d)$, we define
\begin{align*}
V\coloneqq \{ v\in H^1(\Omega;\R^d);~
v|_{\Gamma_D}=0\},
\quad V(b)\coloneqq V+b=\{ v\in H^1(\Omega;\R^d);~
v-b\in V\}.
\end{align*}
Then, setting $g^k\coloneqq g(\cdot,t_k)\in H^1(\Omega;\R^d)$, we will consider that the weak solution $u^k$ for $k\in \{2,\cdots,m\}$ belongs to the affine space $V(g^k)$.

\subsection{Time discretization for DF-PFM}
Under the initial conditions \eqref{uz0}, we consider the following time discretization for the dynamic F-PFM \eqref{DF}:
\begin{subequations}
\begin{empheq}
[left=\empheqlbrace]{align}
& \hspace{5pt} 
\rho\frac{u^k-2 u^{k-1}+u^{k-2}}{\tau^2}=\DIV \sigma [u^k,z^{k-1}] 
&&\mbox{in}~\Omega~~(k=2,3,...m),
\label{d-DF1}\\
& \hspace{5pt} 
\alpha\frac{\tilde{z}^k-z^{k-1}}{\tau} = \epsilon \DIV (\gamma_* \nabla \tilde{z}^k ) - \frac{\gamma_*}{\epsilon}\tilde{z}^k + W^k(1-\tilde{z}^k)
&&\mbox{in}~{\Omega}~~ (k=1,2,...m), \label{d-DF2}\\
& \hspace{5pt} 
\sigma[u^k,z^{k-1}]\nu=0
&&\mbox{on}~\Gamma_N
~~(k=2,3,...m),\label{d-DF3}\\
& \hspace{5pt} u^k=g^k
&&\mbox{on}~ {\Gamma_D}~~(k=2,3,...m),\label{d-DF4}\\
& \hspace{5pt} 
\frac{\partial \tilde{z}^k}{\partial \nu}=0 
&&\mbox{on}~{\Gamma}~~(k=1,2,...m),\label{d-DF5}\\
& \hspace{5pt} z^k=\max (\tilde{z}^k, z^{k-1})
&&\mbox{in}~\Omega ~~(k=1,2,...m),\label{d-DF6}
\end{empheq}\label{d-DF}
\end{subequations}
where we set $W^k\coloneqq W(u^k)$.

The weak form of $u^k$ corresponding to \eqref{d-DF1}, \eqref{d-DF3}, and \eqref{d-DF4} for $k=2,\cdots,m$ is given as follows:
\begin{align}\label{uk}
\begin{cases}
~u^k \in V(g^k)\\
~{\DS \int_\Omega (1-z^{k-1})^2\sigma[u^k] : e[v]\,dx + \frac{\rho}{\tau^2}\int_\Omega u^k \cdot v \,dx = \frac{\rho}{\tau^2}\int_\Omega \tilde{u}^{k-1} \cdot v\, dx \quad (v \in V),} 
\end{cases}
\end{align}
where we set $\tilde{u}^{k-1}\coloneqq 2u^{k-1}-u^{k-2}$.
We remark that the weak solution $u^k$ to \eqref{uk} uniquely exists
if $z^{k-1}\in L^\infty(\Omega)$ and $\tilde{u}^{k-1}\in L^2(\Omega)$ from the Lax-Milgram theorem \cite{Ciarlet1978}.

For $z^k$, the above linear implicit scheme \eqref{d-DF2}, \eqref{d-DF5}, and \eqref{d-DF6} was originally proposed in \cite{TK_2009} as a time-discrete solution method for the irreversible gradient flows. The convergence of the scheme was later proven in \cite{K-N2021} in the case of a simplified irreversible gradient flow. The weak form of $\tilde{z}^k$ for $k=1,\cdots,m$ is given as follows:
\begin{align*}
\begin{cases}
~\tilde{z}^k \in H^1(\Omega)\\
~{\DS \int_\Omega \left(\frac{\alpha}{\tau} + \frac{\gamma_*}{\epsilon} + W^k\right)\tilde{z}^k\zeta\, dx 
+ \epsilon \int_\Omega\gamma_* \nabla \tilde{z}^k\cdot \nabla \zeta \,dx = \int_\Omega \left(\frac{\alpha}{\tau}z^{k-1} + W^k\right)\zeta \,dx \quad (\zeta \in H^1(\Omega)) .} 
\end{cases}
\end{align*}

\subsection{Time discretization for DF-PFM with a unilateral contact condition}
In contrast to the preceding section, it is imperative to develop an effective numerical method for handling the terms $\sigma^\dagger [u,z]$ and $W_+(u)$ in the dynamic F-PFM with a unilateral contact condition \eqref{DFU}. In this paper, we propose the following method. For $u\in H^1(\Omega ; \R^d)$, define the variable $\xi[u]\in L^\infty(\Omega)$ as:
\begin{align*}
\xi[u](x)\coloneqq 
\begin{cases}
1&(\mbox{if}~\DIV u(x)\ge 0),\\
0&(\mbox{if}~\DIV u(x)< 0).
\end{cases}
\end{align*}
The terms $(\DIV u)_\pm$ can be written in the following forms:
\begin{align*}
(\DIV u)_+ =\xi[u]\DIV u,\quad (\DIV u)_- =(1-\xi[u])\DIV u
\quad \mbox{a.e. in }\Omega.
\end{align*}
For $u\in H^1(\Omega;\R^d)$ and $\xi$, $z\in L^\infty (\Omega)$, we define
\begin{align*}
\tilde{\sigma}^\dagger [u,\xi,z] 
&\coloneqq (1-z)^2\big(\,\lambda^* \xi(\DIV u)I+2\mu e_D[u]\,\big)
+\lambda^*(1-\xi)(\DIV u)I,\\
\tilde{W}_+(u,\xi) 
&\coloneqq \lambda^*\xi(\DIV u)^2+\mu |e_{D}[u]|^2.
\end{align*}
Then, we have
\begin{align*}
\tilde{\sigma}^\dagger [u,\xi[u],z]=\sigma^\dagger [u,z],\quad
\tilde{W}_+(u,\xi[u])=W_+(u)~~\mbox{a.e. in }\Omega,
\end{align*}
for $u\in H^1(\Omega;\R^d)$ and $z\in L^\infty (\Omega)$.
Moreover, 
since 
\begin{align*}
\tilde{\sigma}^\dagger [u,\xi,z] 
= \eta (\DIV u)I+2\mu (1-z)^2e[u],
\quad
\eta \coloneqq (1-z)^2\left( \lambda^*\xi-\frac{2\mu}{d}\right)+\lambda^*(1-\xi),
\end{align*}
we have
\begin{align*}
\int_\Omega &\left(\DIV \tilde{\sigma}^\dagger [u,\xi,z] \right)\cdot v\,dx
= \int_{\Gamma_N}\left(\tilde{\sigma}^\dagger [u,\xi,z]\nu\right)\cdot v\,ds
-\int_\Omega \tilde{\sigma}^\dagger [u,\xi,z] : e[v]\,dx\\
&= \int_{\Gamma_N}\left(\tilde{\sigma}^\dagger [u,\xi,z]\nu\right)\cdot v\,ds
-\int_\Omega \left(\eta (\DIV u)(\DIV v)+2\mu (1-z)^2e[u]:e[v]\right)\,dx.\\
\end{align*}

Therefore, under the initial conditions \eqref{uz0}, we consider the following time discretization for the dynamic F-PFM with a unilateral condition \eqref{DFU}:
\begin{subequations}
\begin{empheq}
[left=\empheqlbrace]{align}
& \hspace{5pt} 
\rho\frac{u^k-2 u^{k-1}+u^{k-2}}{\tau^2}=\DIV \tilde{\sigma}^\dagger [u^k,\xi^{k-1},z^{k-1}] 
&&\mbox{in}~\Omega~~(k=2,3,...m),
\label{d-DFU1} \\
& \hspace{5pt} 
\alpha\frac{\tilde{z}^k-z^{k-1}}{\tau} = \epsilon \DIV (\gamma_* \nabla \tilde{z}^k ) - \frac{\gamma_*}{\epsilon}\tilde{z}^k + W_+^k(1-\tilde{z}^k) 
&&\mbox{in}~{\Omega}~~ (k=1,2,...m), \label{d-DFU2}\\
& \hspace{5pt} 
\tilde{\sigma}^\dagger[u^k,\xi^{k-1},z^{k-1}]\nu=0
&&\mbox{on}~\Gamma_N
~~(k=2,3,...m),\label{d-DFU3}\\
& \hspace{5pt} u^k=g^k
&&\mbox{on}~ {\Gamma_D}~~(k=2,3,...m),\label{d-DFU4}\\
& \hspace{5pt} 
\frac{\partial \tilde{z}^k}{\partial \nu}=0 
&&\mbox{on}~{\Gamma}~~(k=1,2,...m),\label{d-DFU5}\\
& \hspace{5pt} z^k=\max (\tilde{z}^k, z^{k-1})
&&\mbox{in}~\Omega ~~(k=1,2,...m),\label{d-DFU6}
\end{empheq}\label{d-DFU} 
\end{subequations}
where we set $\xi^k\coloneqq \xi [u^k]$ and $W_+^k\coloneqq \tilde{W}_+(u^k,\xi^{k-1})$.

The weak form of $u^k$ corresponding to \eqref{d-DFU1}, \eqref{d-DFU3}, and \eqref{d-DFU4} for $k=2,\cdots,m$ is given as follows:
\begin{align}\label{uk2}
\begin{cases}
~u^k \in V(g^k)\\
~{\DS \int_\Omega \left( \eta^{k-1}(\DIV u^k)(\DIV v)+e[u^k] : e[v]\right)\,dx + \frac{\rho}{\tau^2}\int_\Omega u^k v \,dx = \frac{\rho}{\tau^2}\int_\Omega \tilde{u}^{k-1}v\,dx \quad (v \in V),} 
\end{cases}
\end{align}
where we set 
\begin{align*}
\eta^k\coloneqq (1-z^k)^2\left( \lambda^*\xi^k-\frac{2\mu}{d}\right)+\lambda^*(1-\xi^k).
\end{align*}
We remark that \eqref{uk2} is a linear problem and its solution $u^k$ is again uniquely determined by the Lax-Milgram theorem.

For $z^k$, we can just replace $W^k$ by $W_+^k$ in the linear implicit scheme for $z^k$ described in \eqref{d-DF}.

\section{Finite element simulations}
In this section, we present a comparative analysis of two mathematical models: the dynamic fracture phase-field model (DF-PFM) \eqref{DF} and the DF-PFM with a unilateral contact condition \eqref{DFU}. This analysis is carried out through numerical experiments. The simulations were performed using the finite element software FreeFEM \cite{Hecht_2012}, applying the P1 finite element method at each time step to the time-discretized schemes \eqref{d-DF} and \eqref{d-DFU} described in the previous section.

\subsection{Problem setup}
We apply the aforementioned models to investigate crack propagation caused by compressive loading and the introduction of a primary wave within the domain. We begin by describing the computational domain and the initial crack configuration, followed by the numerical results for the damage variable and the dynamic propagation of the crack.

\begin{figure}[!ht]
    \centering
     \begin{subfigure}{0.43\textwidth}
    \centering
    \includegraphics[width=0.8\textwidth]{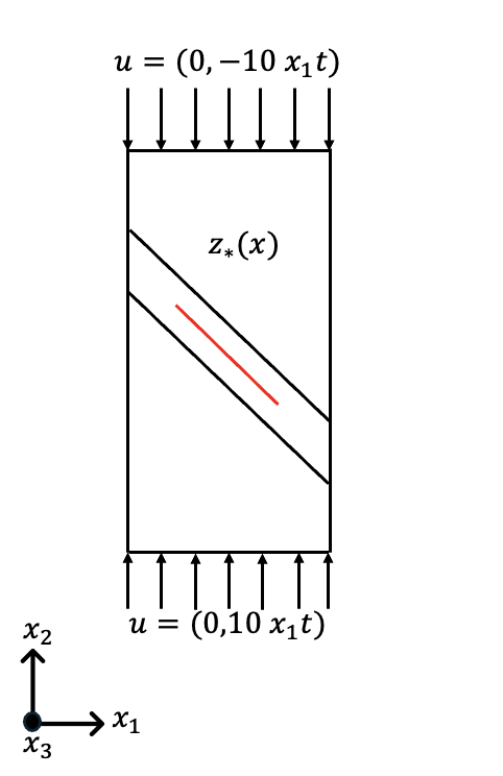}
    \caption{}\label{fig:D1}
    \end{subfigure}
    \begin{subfigure}{0.43\textwidth}
    \centering
    \includegraphics[width=0.7\textwidth]{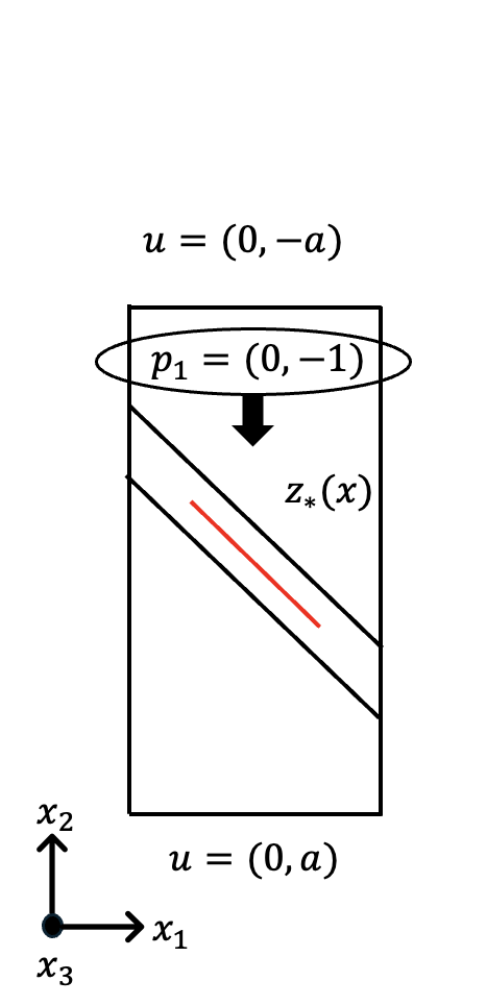}
    \caption{}\label{fig:D2}
    \end{subfigure}
    \caption{Setup of the domain, boundary conditions, and initial crack. (a) Initial compression test (left), (b) P-wave injection under fixed compressive boundary conditions (right).}\label{2DED}
   \end{figure}

Assuming $d=2$, we define the domain as $\Omega=(-1/2,1/2)\times(-1,1)\subset  \R^2$, with boundary conditions specified as follows:
	\begin{align}
    &\Gamma_{D} = \Gamma_{+ D}\cup \Gamma_{-D}, \quad \Gamma_{\pm D} \coloneqq \Gamma \cap \{x_{2} =\pm 1\}, \quad \Gamma_{N}=\Gamma \setminus \Gamma_D.  \notag
	\end{align}
The initial crack is described by the function:
\begin{align}
    z_0(x) \coloneqq \frac{e^{-\eta_1^2/\epsilon^2}}{1+e^{(\eta_2-0.15)/\epsilon}+e^{-(\eta_2+0.15)/\epsilon}}. \notag
\end{align}
where, $\eta_1=x_1\sin \theta+x_2\cos \theta$, $\eta_2=x_1\cos \theta-x_2\sin \theta$ and for this case we consider $\theta=\frac{\pi}{4}$. The Dirichlet boundary condition is given by 
$g(x,t)=(0,-10x_2t)\in \R^2$ for $x=(x_1,x_2)\in \Gamma_D$,
while the Neumann boundary condition is homogeneous: $q=0$ on $\Gamma_N$. Details of the domain and boundary configuration are illustrated in Figures~\ref{fig:D1} and \ref{fig:D2}.  

For the numerical computations, we use the following nondimensional parameter values: Young’s modulus ($E_Y$) is set to $50$, Poisson’s ratio ($\nu_P$) to $0.29$, fracture energy ($\gamma_*$) to $0.5$, density of the elastic body ($\rho$) to $5.0 \times 10^{-4}$, time relaxation constant ($\alpha$) to $10^{-4}$, time increment ($\tau$) to $2 \times 10^{-5}$, and internal length scale ($\epsilon$) to $0.01$.

We first apply a time-dependent displacement to induce compression, as shown in the domain in Fig.~\ref{fig:D1}, and observe the displacement value immediately before crack initiation using the dynamic elastic wave equation model.

Next, we perform a numerical experiment in which a primary wave (P-wave) is introduced into the domain to impact the crack. This is done using both the standard DF-PFM and the DF-PFM with a unilateral contact condition. For this simulation, we use the displacement field immediately before crack initiation from the compression test as the initial condition, fixing it to a value of $a = 0.240$, as shown in Fig.~\ref{fig:D2}. The P-wave is injected from top to bottom, i.e., in the direction $p_1 = (0, -1)$.
Setting the speed of the P-wave $v_p\coloneqq \sqrt{(\lambda+2 \mu)/\rho}$,
we define
\begin{align*}
u_\text{ini}(x,t)\coloneqq 
\left(0,~0.01*\exp \left(-\left(\frac{x_2-v_pt-0.5}{0.1}\right)^2\right)\right)^T.
\end{align*}
This is a plane P-wave in the $p_1$ direction.
for the case of no cracks or damage.
In the numerical example, we set
\begin{align*}
u^0(x):=u_\text{ini}(x,0),\quad
u^1(x):=u_\text{ini}(x,\tau).
\end{align*}
In this case, the boundary conditions and initial conditions are not consistent, but we will assume that the compression boundary conditions begin at $k=2$ and thereafter.

\subsection{Wave propagation analysis}
We plot the squared magnitude of the wave speed, $|v|^2$, where the velocity is approximated by $v \coloneqq (u^k - u^{k-1}) / \tau$.
Figures~\ref{2DED1-1}-\ref{2DED1-3} illustrate the wave propagation for both models.

In Figure~\ref{2DED1-1} ($t=0.0\sim 1.2\times 10^{-3}$), the injected P-wave passes near the crack without much effect on the initial crack under compression, which is observed for both DF-PFM (top) and DF-PFM with a unilateral contact condition (bottom) with little difference. On the other hand, in Figure~\ref{2DED1-2} ($t=1.6\times 10^{-3}\sim 3.2\times 10^{-3}$), no significant change is observed in DF-PFM  with a unilateral contact condition, but in DF-PFM, crack propagation starts and new elastic waves are generated due to the singularity of the crack tip. In Figure~\ref{2DED1-3} ($t=3.6\times 10^{-3}\sim 5.2\times 10^{-3}$), it can be observed that new elastic waves originating from the singularity of crack propagation are produced in DF-PFM  with a unilateral contact condition as well.

\begin{figure}[!ht]
    \centering
     \begin{subfigure}{0.24\textwidth}
    \centering
    \includegraphics[width=0.8\textwidth]{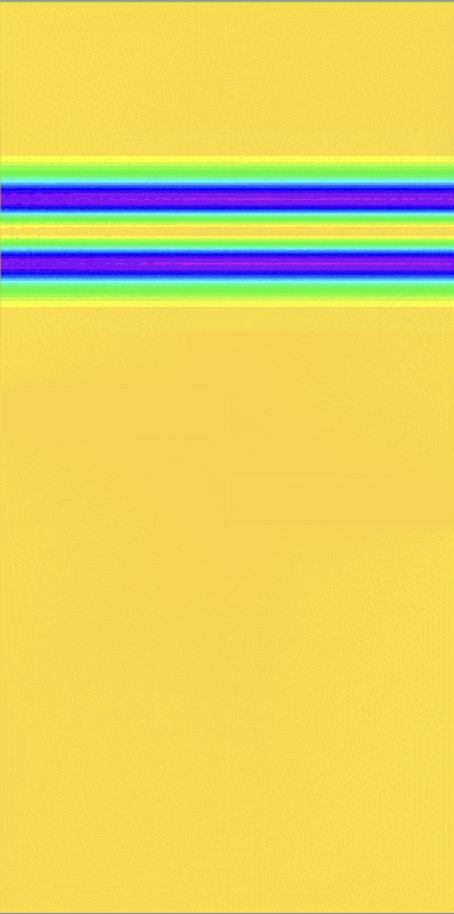}
    \end{subfigure}
    \begin{subfigure}{0.24\textwidth}
    \centering
    \includegraphics[width=0.8\textwidth]{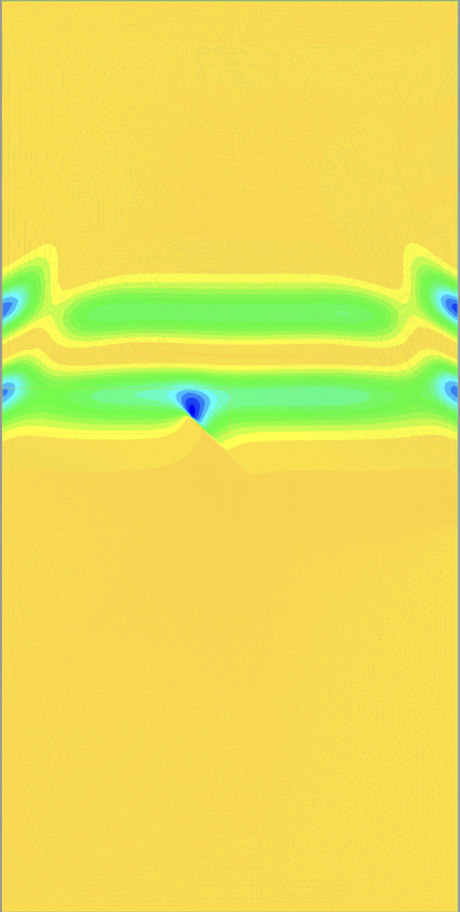}
    \end{subfigure}
    \begin{subfigure}{0.24\textwidth}
    \centering
    \includegraphics[width=0.8\textwidth]{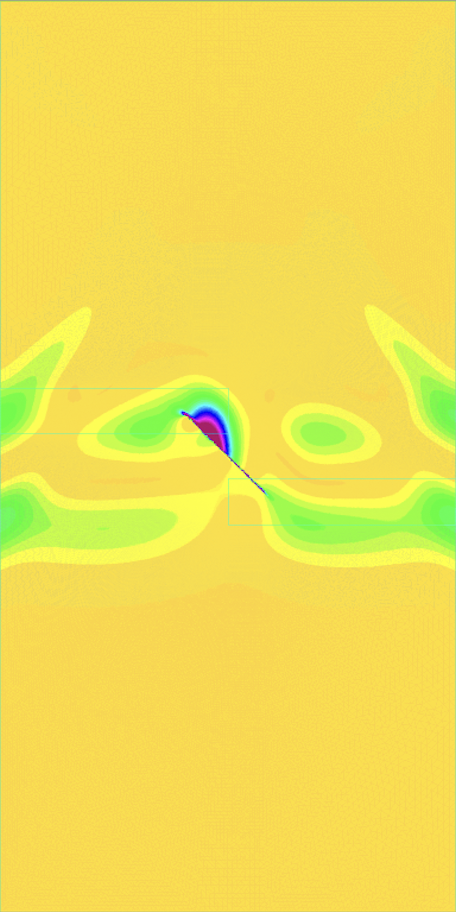}
    \end{subfigure}
    \begin{subfigure}{0.24\textwidth}
    \centering
    \includegraphics[width=0.8\textwidth]{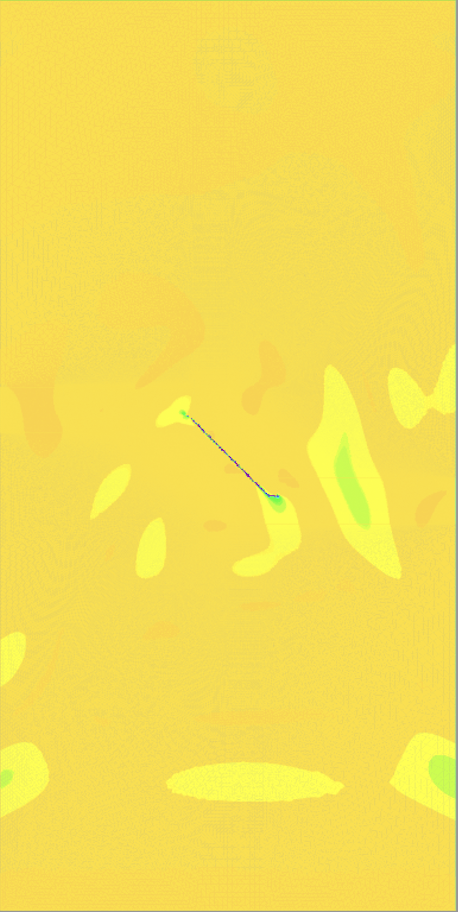}
    \end{subfigure}\\ \vspace{4pt}
  \begin{subfigure}{0.24\textwidth}
    \centering
    \includegraphics[width=0.8\textwidth]{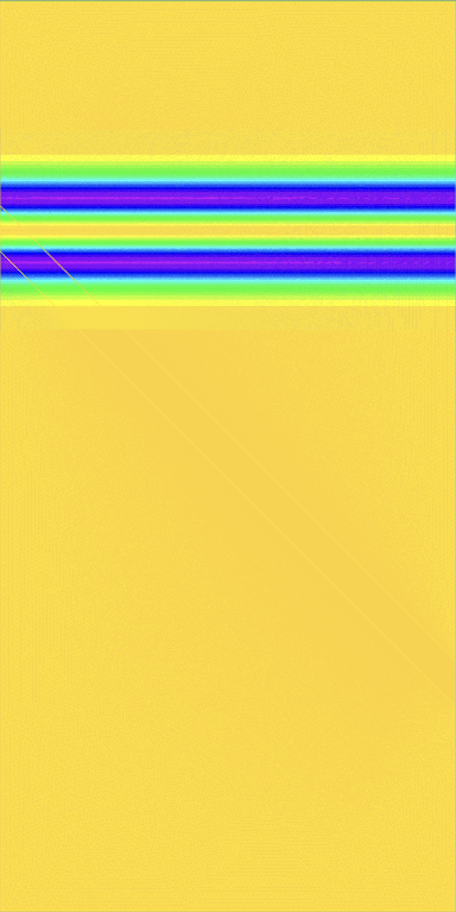}
    \caption{$t=0.0$}
    \end{subfigure}
        \begin{subfigure}{0.24\textwidth}
    \centering
    \includegraphics[width=0.8\textwidth]{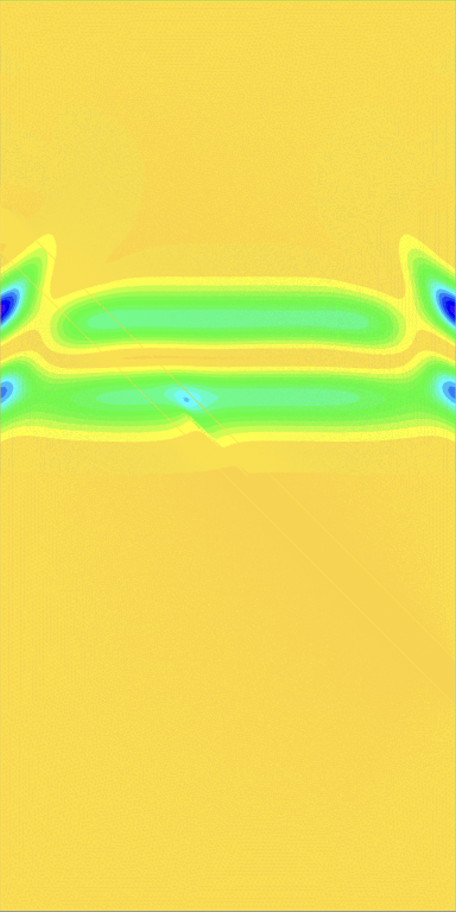}
    \caption{$t=4.0\times10^{-4}$}
    \end{subfigure}
        \begin{subfigure}{0.24\textwidth}
    \centering
    \includegraphics[width=0.8\textwidth]{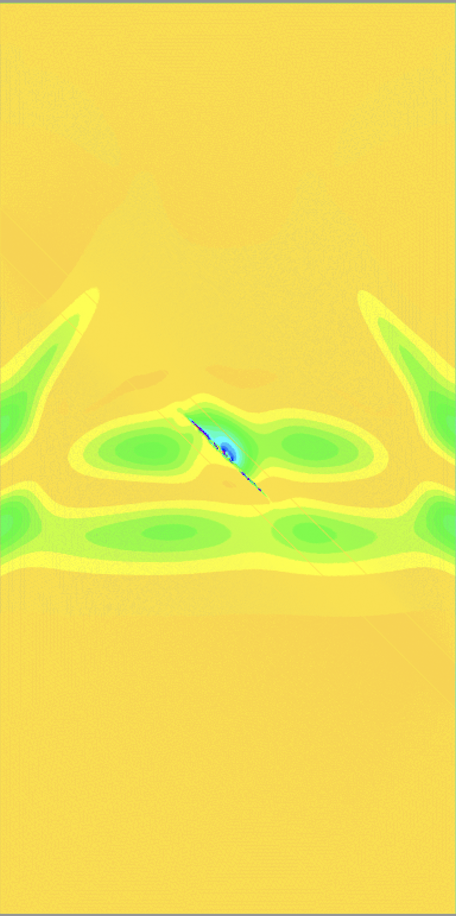}
    \caption{$t=8.0\times10^{-4}$}
    \end{subfigure}
    \begin{subfigure}{0.24\textwidth}
    \centering
    \includegraphics[width=0.8\textwidth]{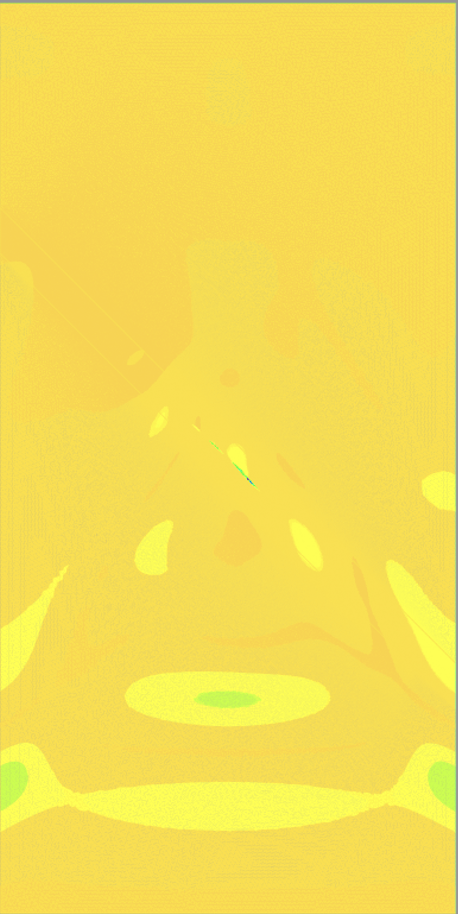}
    \caption{$t=1.2\times10^{-3}$}
    \end{subfigure} \\ \vspace{4pt}
    \begin{subfigure}{0.5\textwidth}
    \centering
    \includegraphics[width=0.8\textwidth]{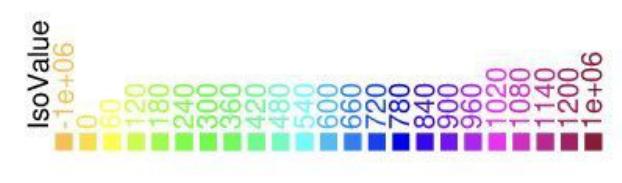}
    \end{subfigure}
    \caption{Simulation results of wave propagation using DF-PFM (top) and DF-PFM with a unilateral contact condition (bottom), with P-wave injection and compression applied at the Dirichlet boundary for $t=0.0$, 
    $4.0\times 10^{-4}$, $8.0\times 10^{-4}$, $1.2\times 10^{-3}$. The color plot represents the value of $|v|^2$.} \label{2DED1-1}
   \end{figure} 

\begin{figure}[!ht]
    \centering
     \begin{subfigure}{0.24\textwidth}
    \centering
    \includegraphics[width=0.8\textwidth]{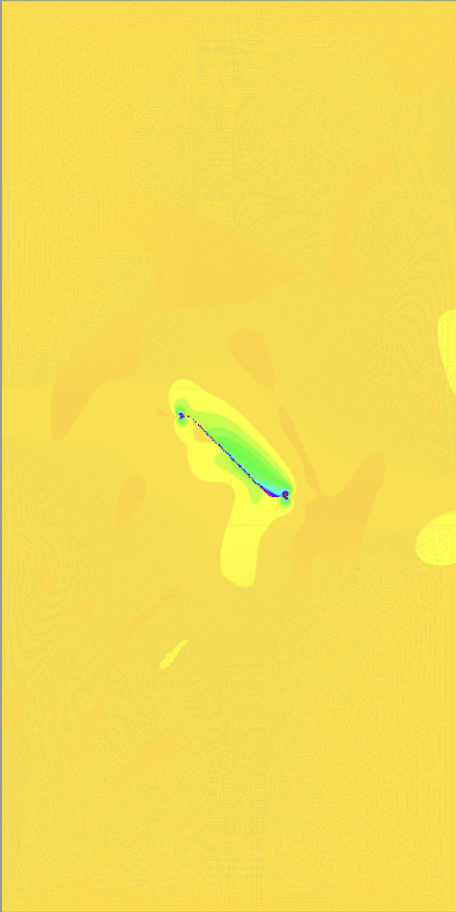}
    \end{subfigure}
    \begin{subfigure}{0.24\textwidth}
    \centering
    \includegraphics[width=0.8\textwidth]{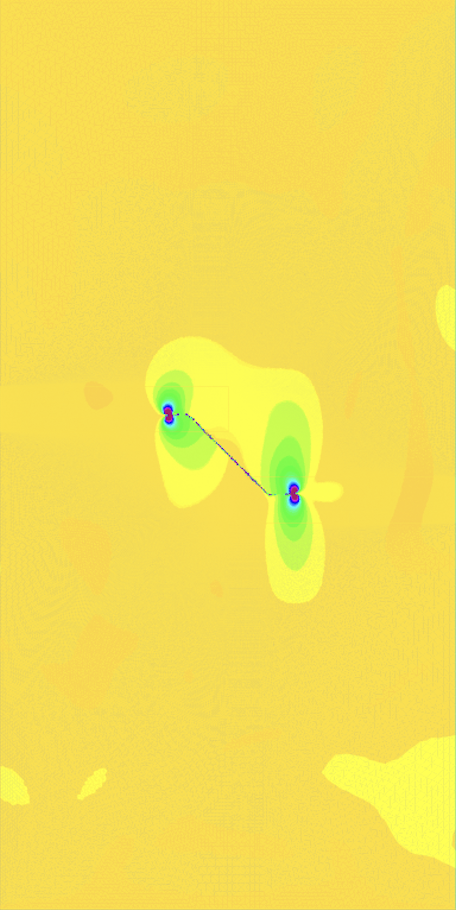}
    \end{subfigure}
    \begin{subfigure}{0.24\textwidth}
    \centering
    \includegraphics[width=0.8\textwidth]{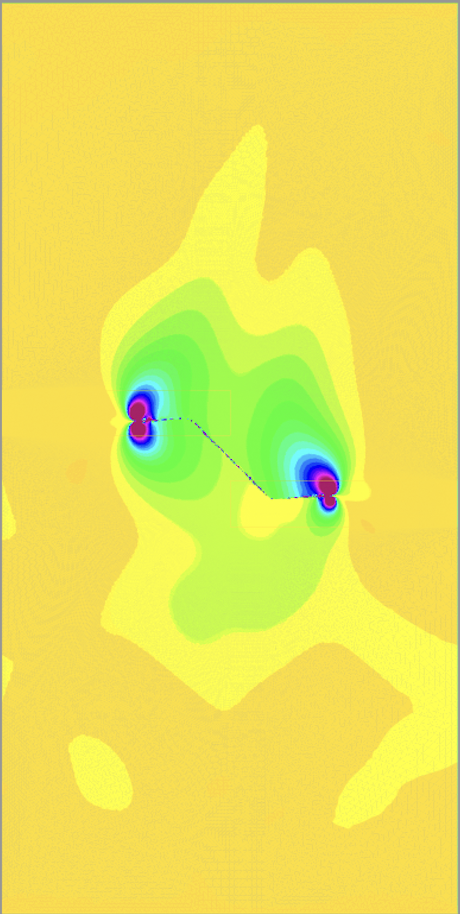}
    \end{subfigure}
    \begin{subfigure}{0.24\textwidth}
    \centering
    \includegraphics[width=0.8\textwidth]{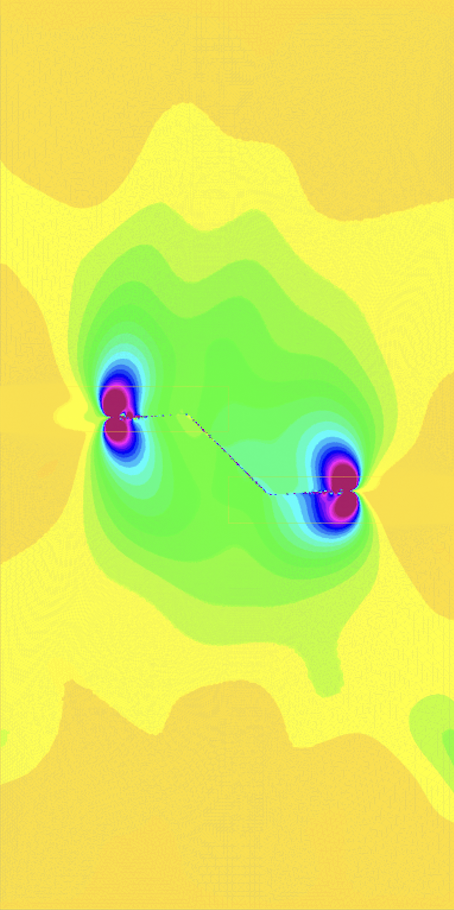}
    \end{subfigure}\\ \vspace{4pt}
  \begin{subfigure}{0.24\textwidth}
    \centering
    \includegraphics[width=0.8\textwidth]{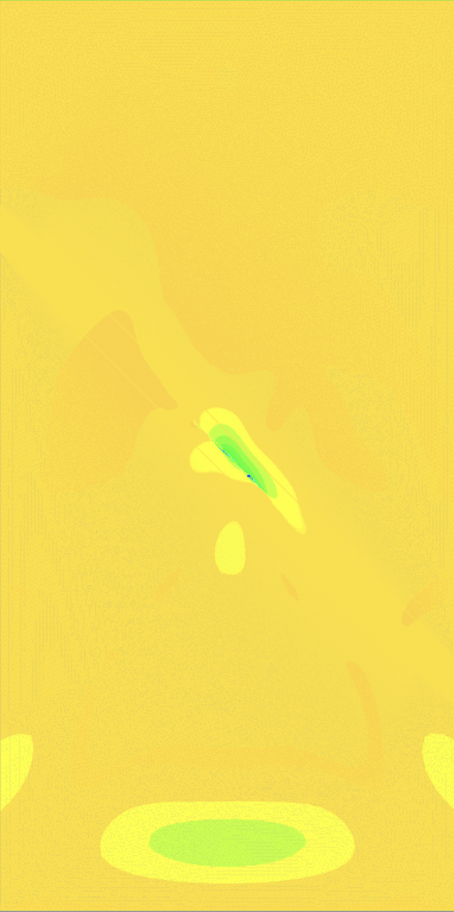}
    \caption{$t=1.6\times10^{-3}$}
    \end{subfigure}
        \begin{subfigure}{0.24\textwidth}
    \centering
    \includegraphics[width=0.8\textwidth]{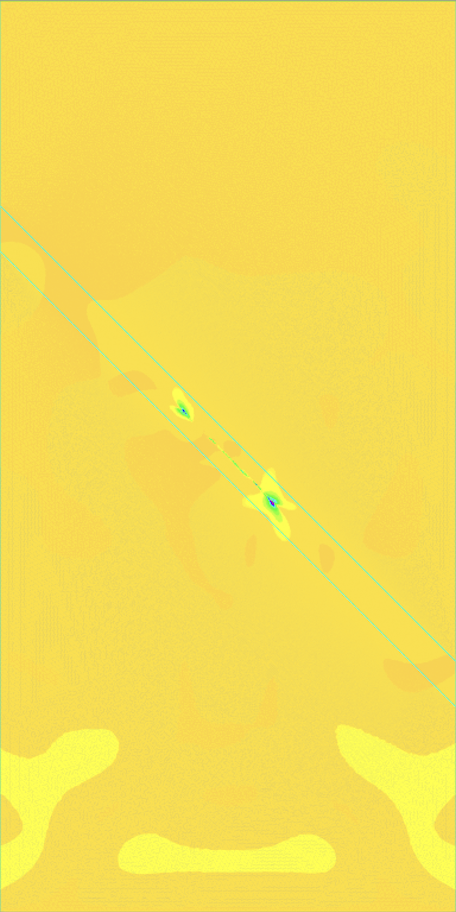}
    \caption{$t=2.0\times10^{-3}$}
    \end{subfigure}
        \begin{subfigure}{0.24\textwidth}
    \centering
    \includegraphics[width=0.8\textwidth]{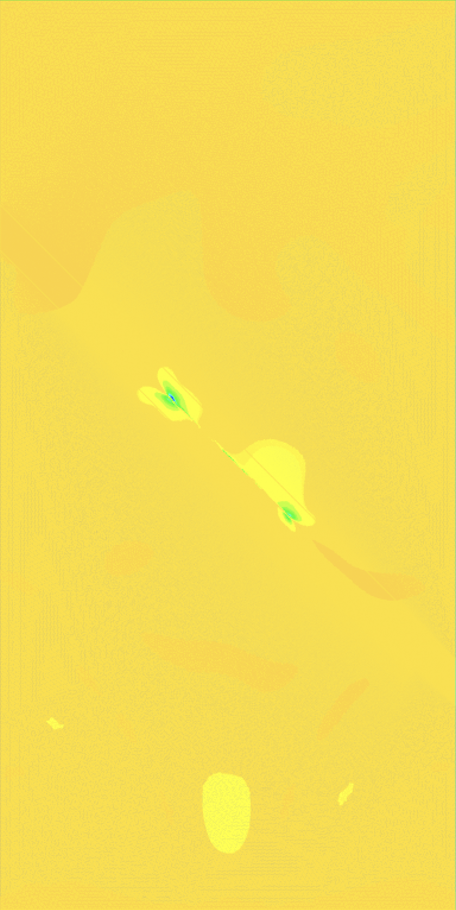}
    \caption{$t=2.8\times10^{-3}$}
    \end{subfigure}
    \begin{subfigure}{0.24\textwidth}
    \centering
    \includegraphics[width=0.8\textwidth]{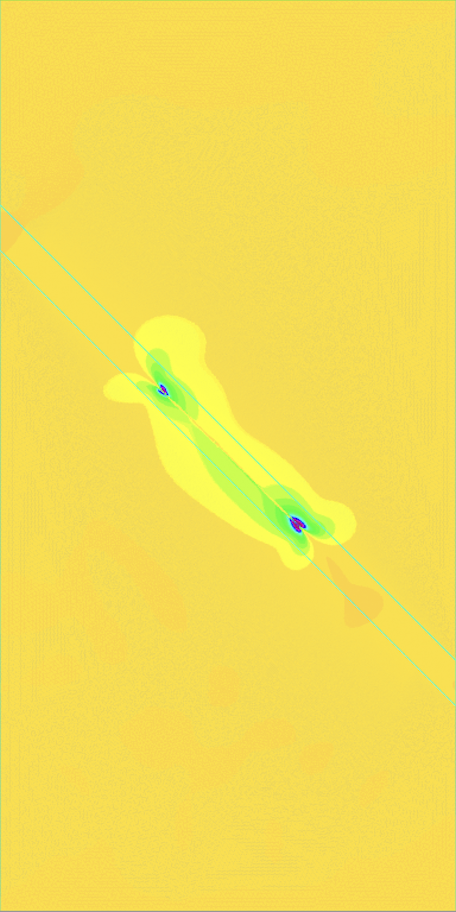}
    \caption{$t=3.2\times10^{-3}$}
    \end{subfigure} \\ \vspace{4pt}
    \begin{subfigure}{0.5\textwidth}
    \centering
    \includegraphics[width=0.8\textwidth]{./fig/ColorMap1.pdf}
    \end{subfigure}
    \caption{(Cont.) Simulation results of wave propagation using DF-PFM (top) and DF-PFM with a unilateral contact condition (bottom), with P-wave injection and compression applied at the Dirichlet boundary for $t=1.6\times 10^{-3}$, $2.0\times 10^{-3}$, $2.8\times 10^{-3}$, $3.2\times 10^{-3}$.The color plot represents the value of $|v|^2$.}\label{2DED1-2}
   \end{figure} 

\begin{figure}[!ht]
    \centering
     \begin{subfigure}{0.24\textwidth}
    \centering
    \includegraphics[width=0.8\textwidth]{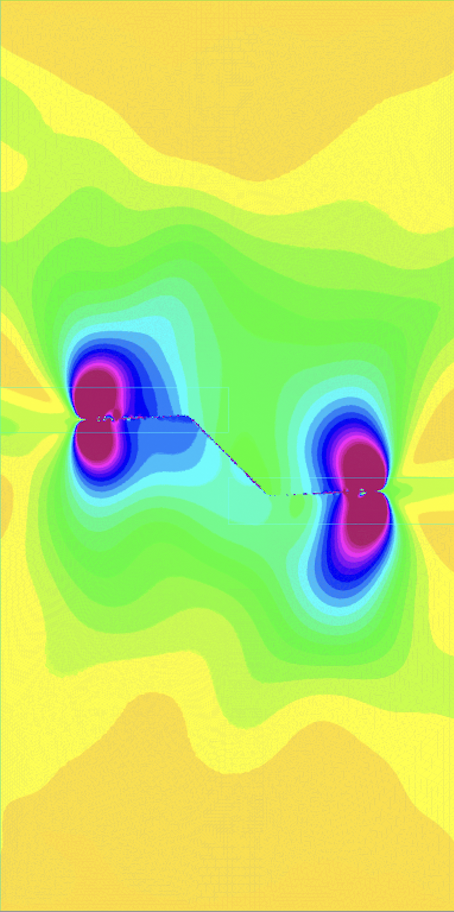}
    \end{subfigure}
    \begin{subfigure}{0.24\textwidth}
    \centering
    \includegraphics[width=0.8\textwidth]{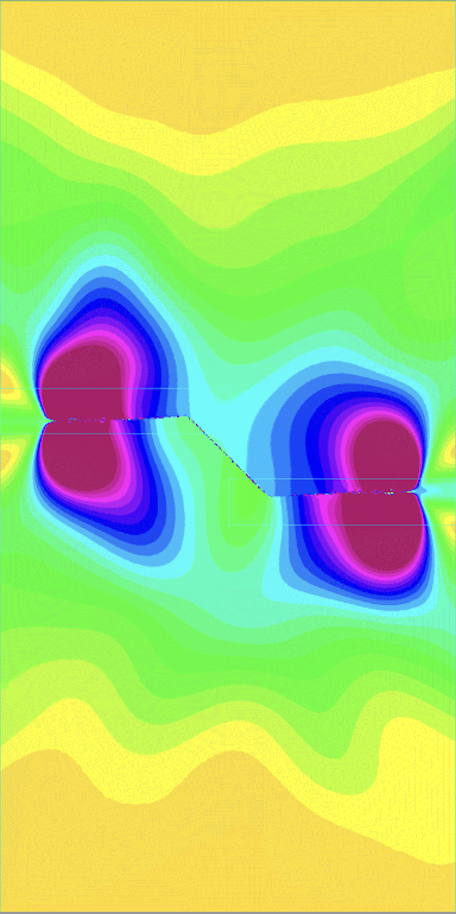}
    \end{subfigure}
    \begin{subfigure}{0.24\textwidth}
    \centering
    \includegraphics[width=0.8\textwidth]{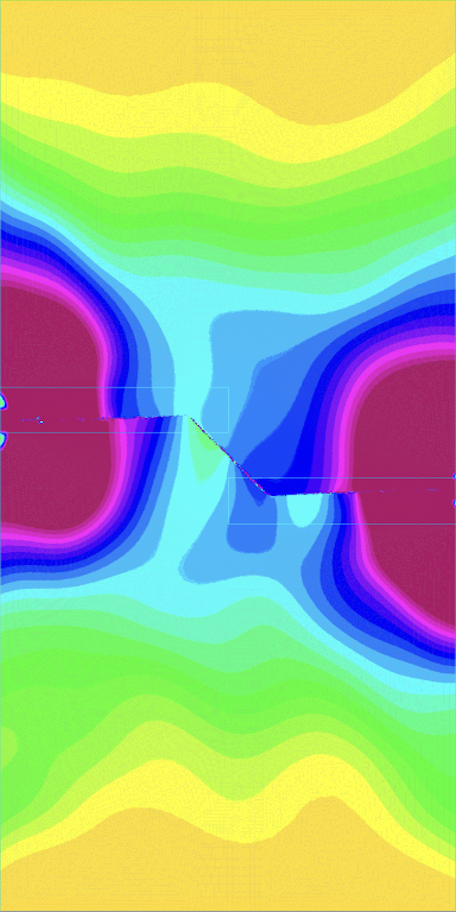}
    \end{subfigure}
    \begin{subfigure}{0.24\textwidth}
    \centering
    \includegraphics[width=0.8\textwidth]{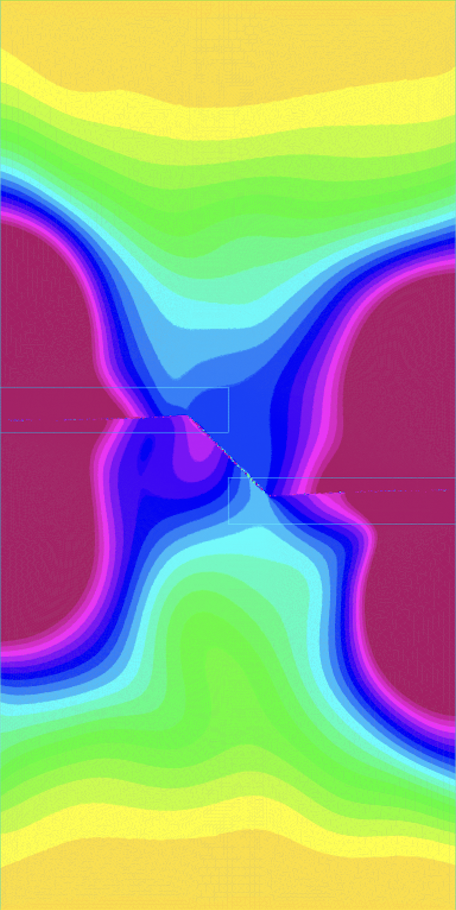}
    \end{subfigure}\\ \vspace{4pt}
  \begin{subfigure}{0.24\textwidth}
    \centering
    \includegraphics[width=0.8\textwidth]{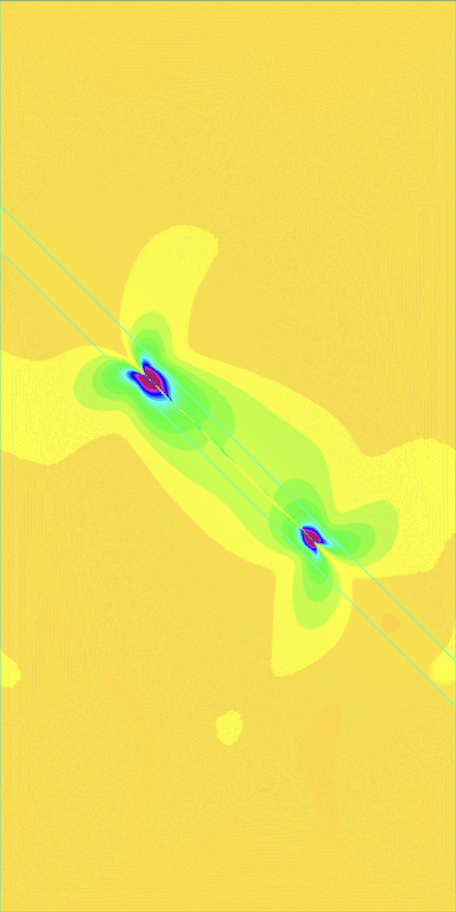}
    \caption{$t=3.6\times10^{-3}$}
    \end{subfigure}
        \begin{subfigure}{0.24\textwidth}
    \centering
    \includegraphics[width=0.8\textwidth]{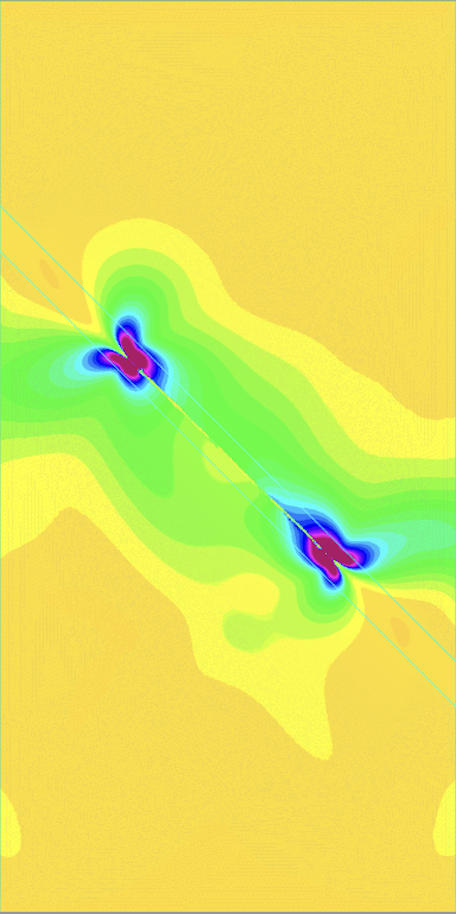}
    \caption{$t=4.4\times10^{-3}$}
    \end{subfigure}
        \begin{subfigure}{0.24\textwidth}
    \centering
    \includegraphics[width=0.8\textwidth]{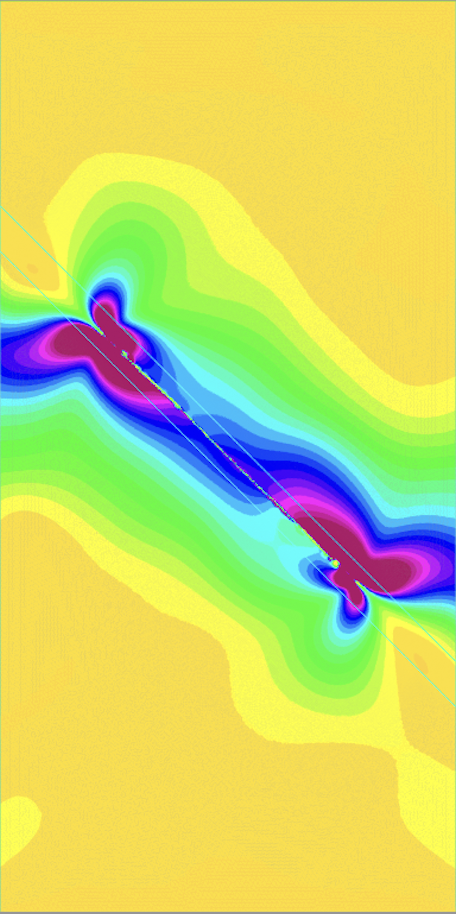}
    \caption{$t=4.8\times10^{-3}$}
    \end{subfigure}
    \begin{subfigure}{0.24\textwidth}
    \centering
    \includegraphics[width=0.8\textwidth]{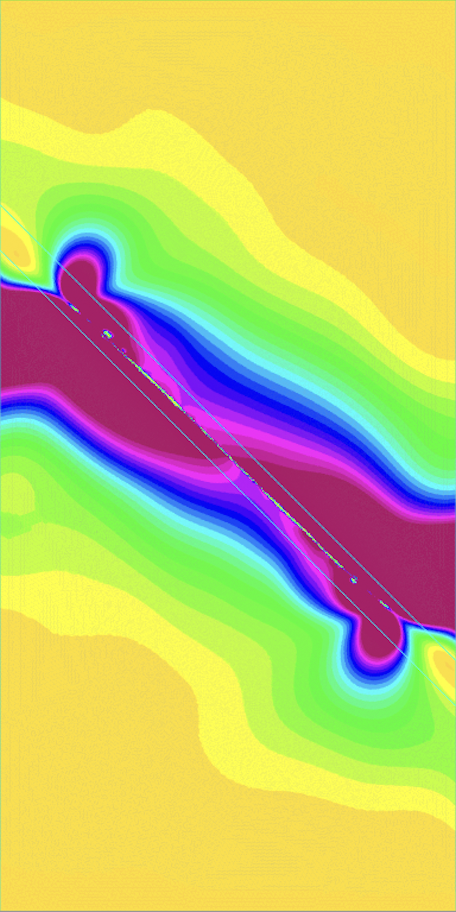}
    \caption{$t=5.2\times10^{-3}$}
    \end{subfigure} \\ \vspace{4pt}
    \begin{subfigure}{0.5\textwidth}
    \centering
    \includegraphics[width=0.8\textwidth]{./fig/ColorMap1.pdf}
    \end{subfigure}
    \caption{(Cont.) Simulation results of wave propagation using DF-PFM (top) and DF-PFM with a unilateral contact condition (bottom), with P-wave injection and compression applied at the Dirichlet boundary for $t=3.6\times 10^{-3}$, $4.4\times 10^{-3}$, $4.8\times 10^{-3}$, $5.2\times 10^{-3}$. The color plot represents the value of $|v|^2$.}\label{2DED1-3}
   \end{figure} 
\subsection{Crack propagation comparison}

In this section, we will examine the simulations from the previous section, focusing on the timing of starting the crack propagation and the morphology of the resulting fractures. Figures~\ref{2DED2-1} and \ref{2DED2-2}, show the resulting crack propagation under compression for each model. 

The crack propagation occurs at approximately $t=1.6 \times 10^{-3}$ for the DF-PFM (top) and around $t=2.8\times 10^{-3}$ for the DF-PFM with a unilateral contact condition (bottom). These observations suggest that the propagation of cracks is delayed in the DF-PFM with a unilateral contact condition relative to the DF-PFM. With regard to the crack morphology, kink-type crack propagation was observed in DF-PFM, with the cracks extending horizontally. Conversely, the oblique direction in DF-PFM with a unilateral contact condition exhibits a fault-slip type crack growth.

\begin{figure}[!ht]
    \centering
     \begin{subfigure}{0.24\textwidth}
    \centering
    \includegraphics[width=0.8\textwidth]{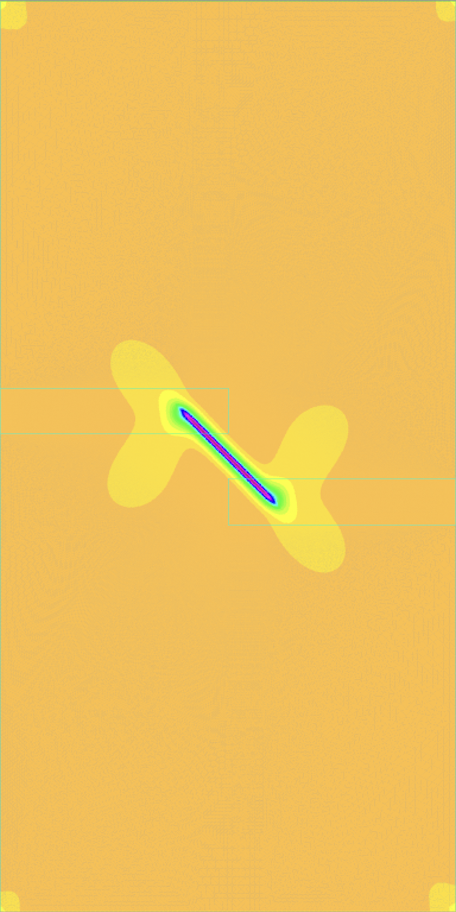}
    \end{subfigure}
    \begin{subfigure}{0.24\textwidth}
    \centering
    \includegraphics[width=0.8\textwidth]{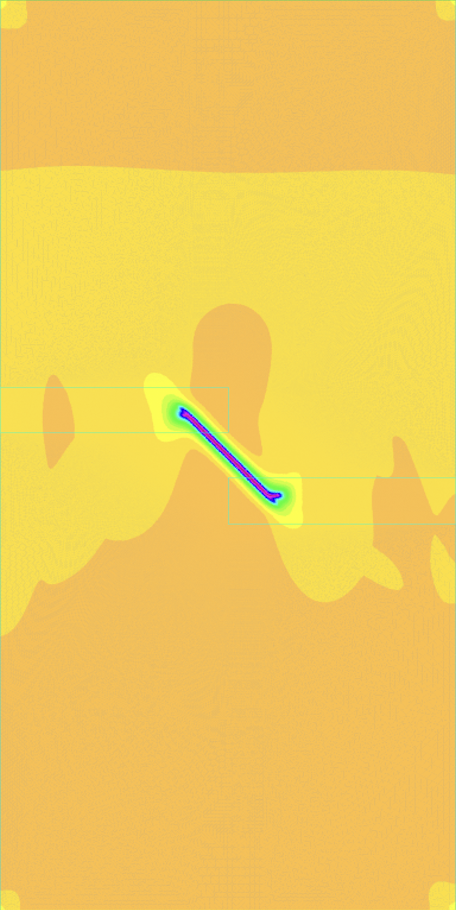}
    \end{subfigure}
    \begin{subfigure}{0.24\textwidth}
    \centering
    \includegraphics[width=0.8\textwidth]{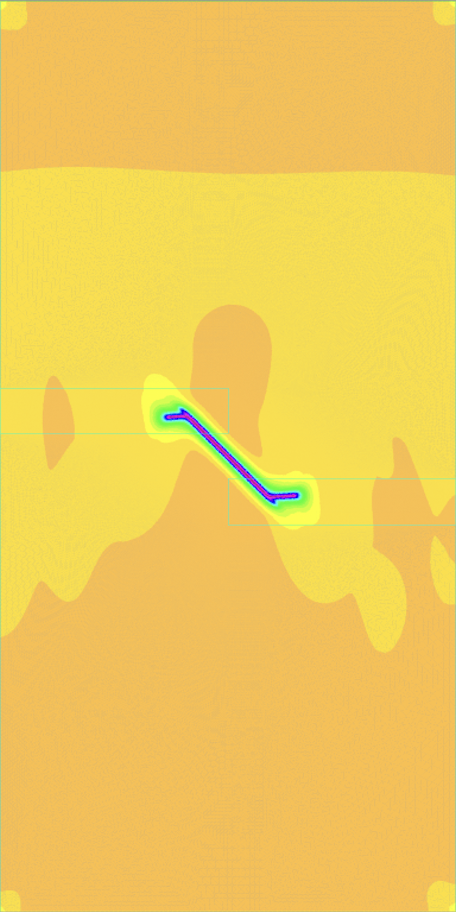}
    \end{subfigure}
    \begin{subfigure}{0.24\textwidth}
    \centering
    \includegraphics[width=0.8\textwidth]{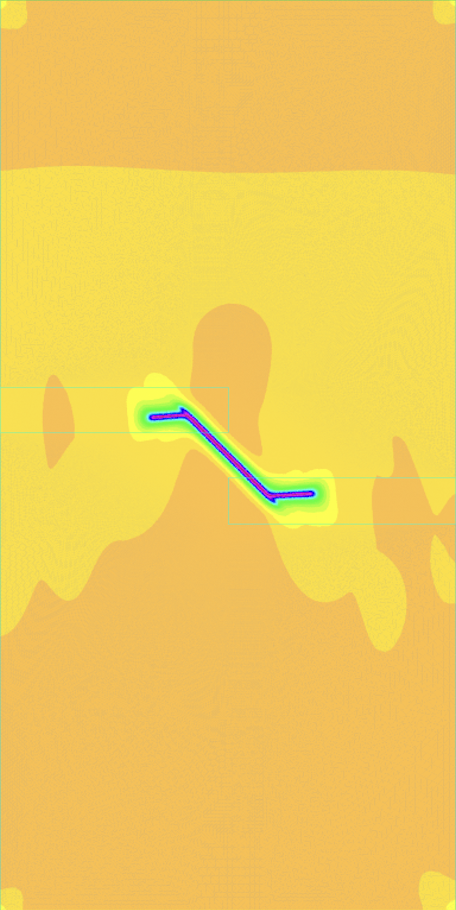}
    \end{subfigure}\\ \vspace{4pt}
  \begin{subfigure}{0.24\textwidth}
    \centering
    \includegraphics[width=0.8\textwidth]{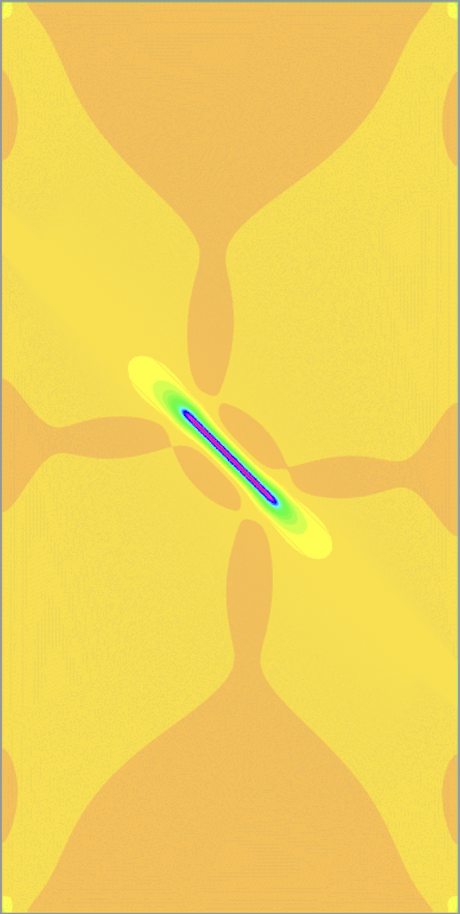}
    \caption{$t=0.0$}
    \end{subfigure}
        \begin{subfigure}{0.24\textwidth}
    \centering
    \includegraphics[width=0.8\textwidth]{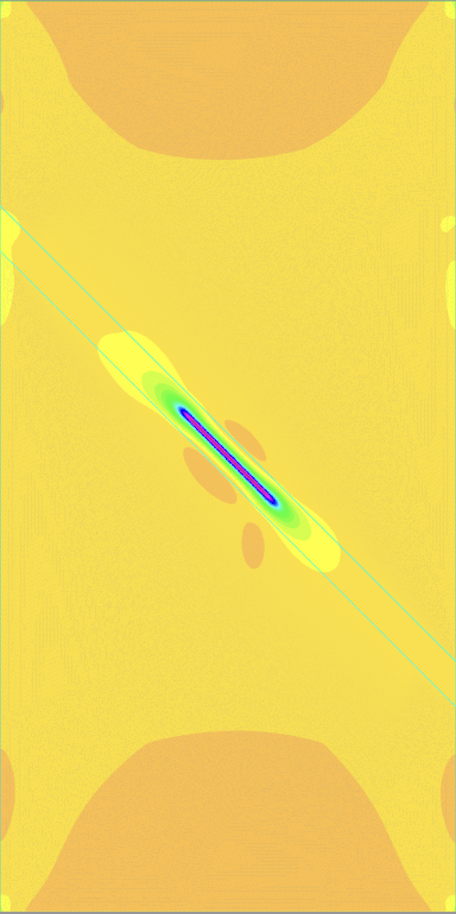}
    \caption{$t=1.6\times10^{-3}$}
    \end{subfigure}
     \begin{subfigure}{0.24\textwidth}
    \centering
    \includegraphics[width=0.8\textwidth]{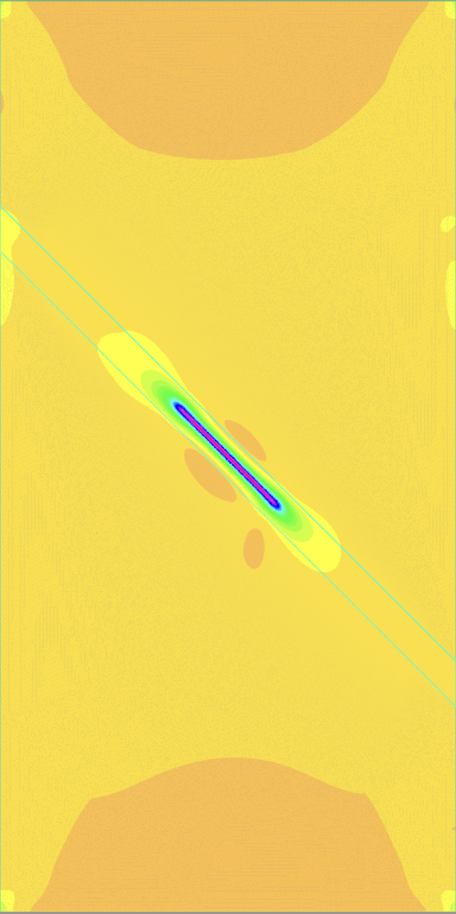}
    \caption{$t=2.4\times10^{-3}$}
    \end{subfigure}
     \begin{subfigure}{0.24\textwidth}
    \centering
    \includegraphics[width=0.8\textwidth]{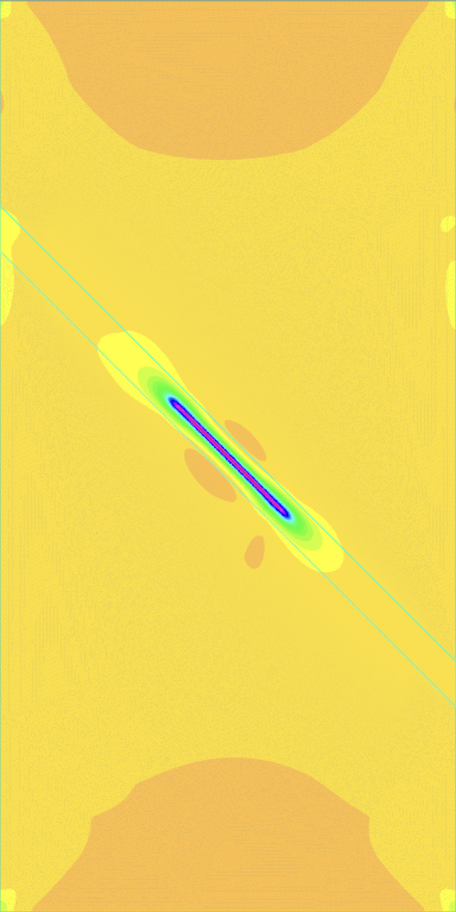}
    \caption{$t=2.8\times10^{-3}$}
    \end{subfigure} \\ \vspace{4pt}
    \begin{subfigure}{0.55\textwidth}
    \centering
    \includegraphics[width=0.8\textwidth]{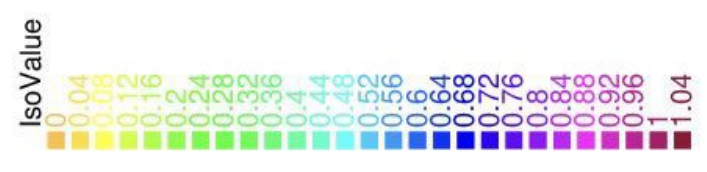}
    \end{subfigure}
    \caption{Simulation results of crack propagation using DF-PFM (top) and DF-PFM with a unilateral contact condition (bottom), under P-wave injection and compression applied at the Dirichlet boundary for $t=0.0$, 
    $1.6\times 10^{-3}$, $2.4\times 10^{-3}$, $2.8\times 10^{-3}$. The color plot represents the value of $z$.}
    \label{2DED2-1}
   \end{figure} 

\begin{figure}[!ht]
    \centering
     \begin{subfigure}{0.24\textwidth}
    \centering
    \includegraphics[width=0.8\textwidth]{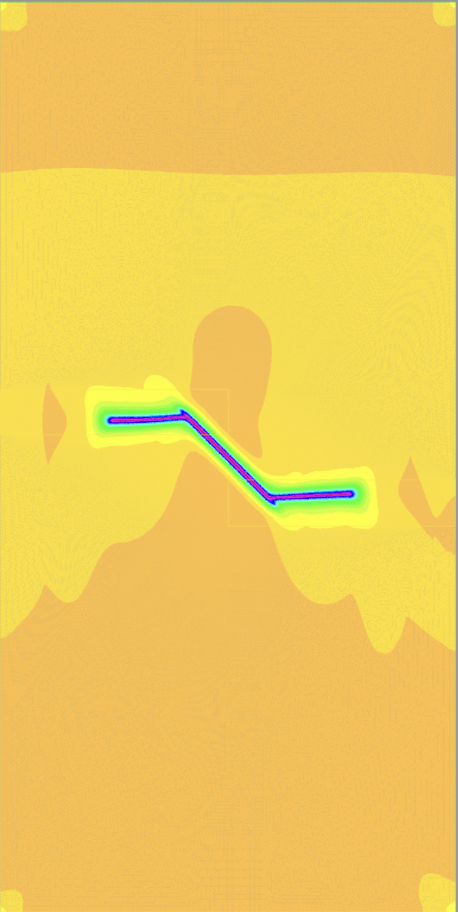}
    \end{subfigure}
    \begin{subfigure}{0.24\textwidth}
    \centering
    \includegraphics[width=0.8\textwidth]{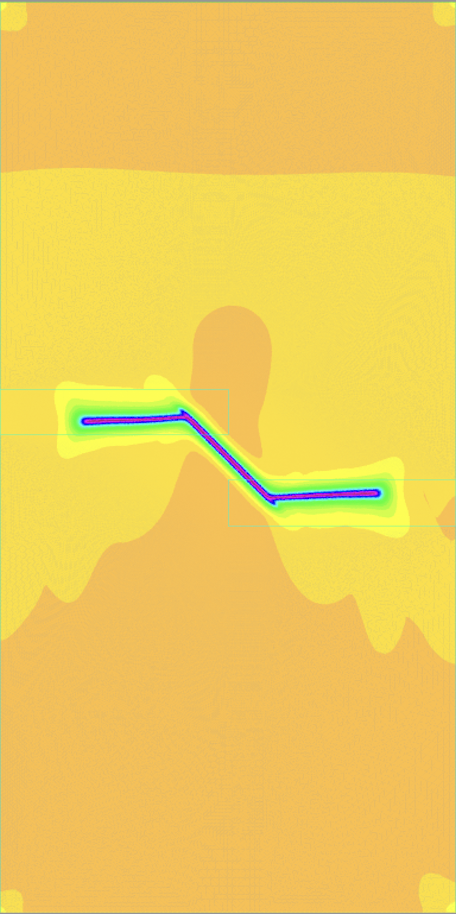}
    \end{subfigure}
    \begin{subfigure}{0.24\textwidth}
    \centering
    \includegraphics[width=0.8\textwidth]{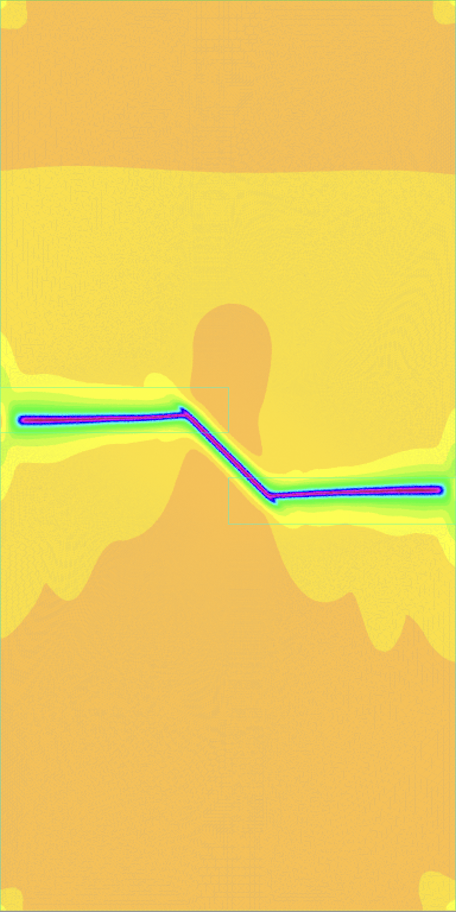}
    \end{subfigure}
    \begin{subfigure}{0.24\textwidth}
    \centering
    \includegraphics[width=0.8\textwidth]{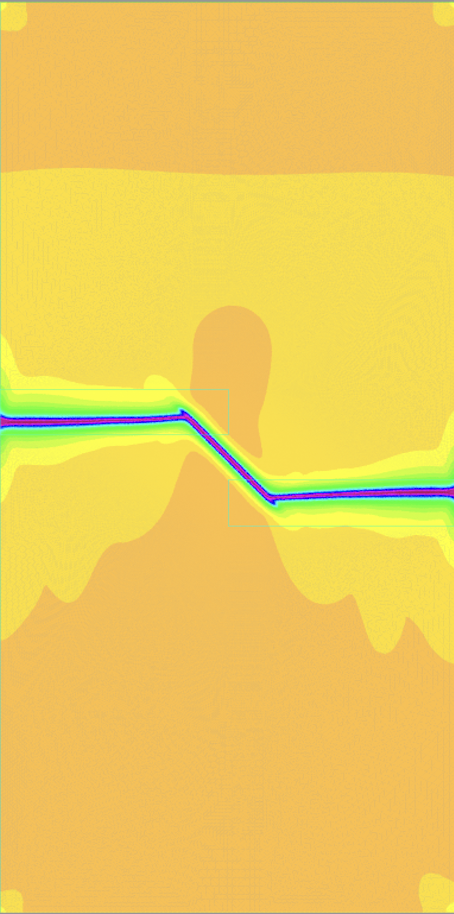}
    \end{subfigure}\\ \vspace{4pt}
  \begin{subfigure}{0.24\textwidth}
    \centering
    \includegraphics[width=0.8\textwidth]{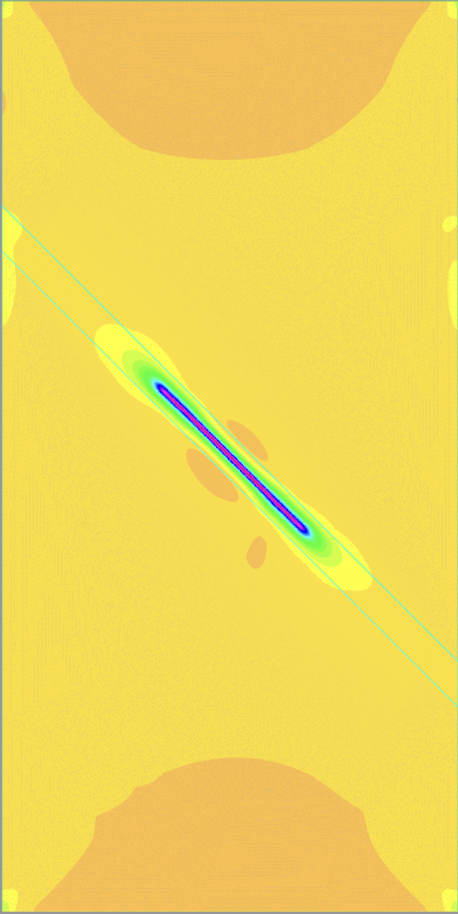}
    \caption{$t=3.6\times10^{-3}$}
    \end{subfigure}
        \begin{subfigure}{0.24\textwidth}
    \centering
    \includegraphics[width=0.8\textwidth]{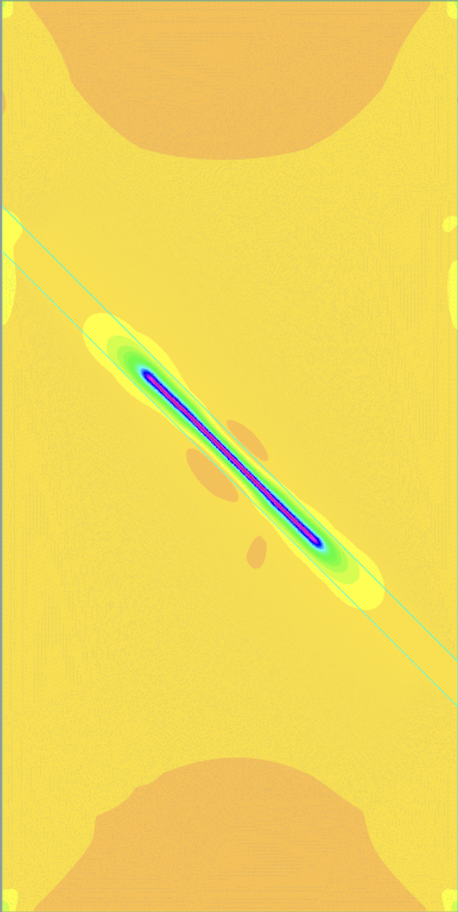}
    \caption{$t=4.0\times10^{-3}$}
    \end{subfigure}
     \begin{subfigure}{0.24\textwidth}
    \centering
    \includegraphics[width=0.8\textwidth]{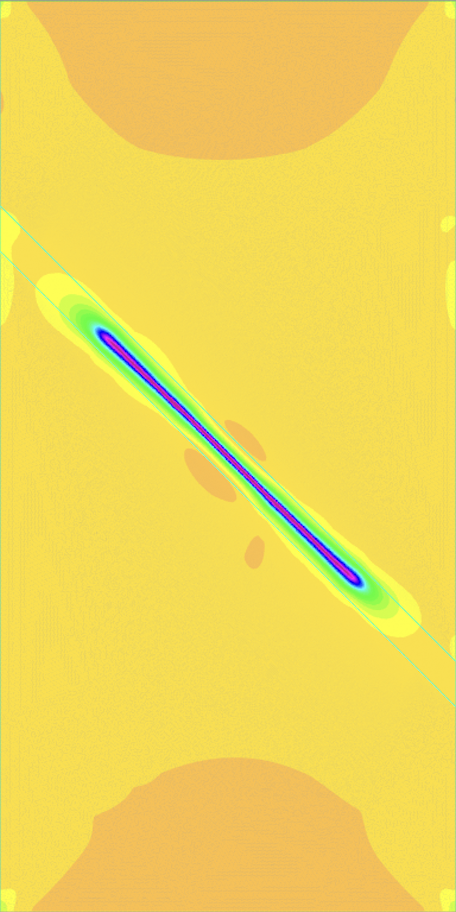}
    \caption{$t=4.8\times10^{-3}$}
    \end{subfigure}
     \begin{subfigure}{0.24\textwidth}
    \centering
    \includegraphics[width=0.8\textwidth]{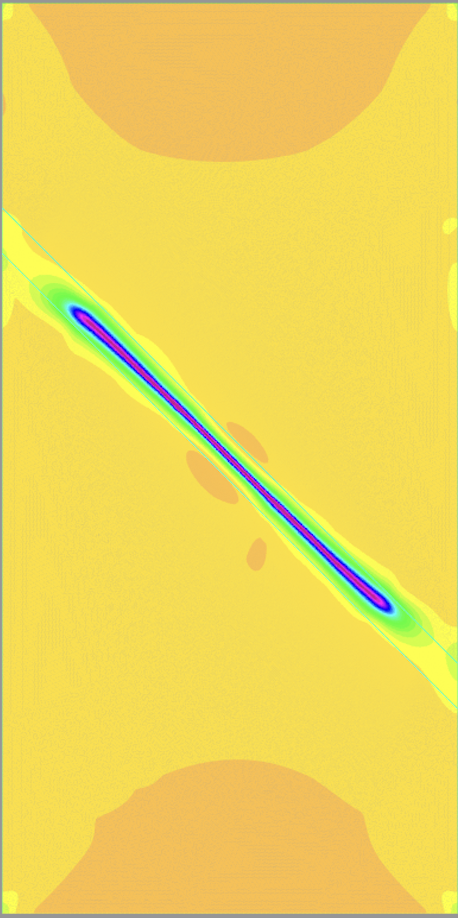}
    \caption{$t=5.2\times10^{-3}$}
    \end{subfigure} \\ \vspace{4pt}
    \begin{subfigure}{0.55\textwidth}
    \centering
    \includegraphics[width=0.8\textwidth]{./fig/ColorMap2.pdf}
    \end{subfigure}
    \caption{(Cont.) Simulation results of crack propagation using DF-PFM (top) and DF-PFM with a unilateral contact condition (bottom), under P-wave injection and compression applied at the Dirichlet boundary for $t=3.6\times 10^{-3}$, $4.0\times 10^{-3}$, $4.8\times 10^{-3}$, $5.2\times 10^{-3}$. The color plot represents the value of $z$.}\label{2DED2-2}
   \end{figure}
\section{Discussion of Results}

In this section, we discuss our results in detail and highlight three key findings: (i) the formulation and theoretical proof of the energy dissipation identity; (ii) the characterization of wave motion in both models; and (iii) the crack propagation behavior in the dynamic fracture phase-field model (DF-PFM) with and without the unilateral contact condition. Each point is elaborated in the following paragraphs.

\paragraph{Energy Dissipation Identity.}
One of the main theoretical contributions of this study is the formal derivation of the energy dissipation identity for both the standard DF-PFM and the modified DF-PFM with a unilateral contact condition. 
This identity ensures that the total energy, composed of elastic energy, kinetic energy, and surface energy, 
together with the energy injection from the external forces, dissipates in line with the second law of thermodynamics. 
Notably, the modified model with the unilateral constraint preserves the energy structure while eliminating non-physical contributions from compressive stresses, which supports its validity as a physically meaningful modification.

\paragraph{Wave Propagation Analysis.}
We applied compressive displacement on the Dirichlet boundary and considered an initially angled crack within the domain. For both models, we recorded the displacement value immediately before the onset of crack propagation. This displacement was then held fixed, and a P-wave was introduced propagating vertically downward. As shown in Figures~\ref{2DED1-1}-\ref{2DED1-3}, the wave behavior was captured numerically for both models. 

In both models, we observe a propagation of the P-wave which passes through the initial crack, and a new elastic wave generated by the singularity of the propagating crack tip. In both cases, it can be inferred that the passage of the P-wave does not seem to affect the crack shape at first glance, but is then triggered by the elastic wave hitting the crack, causing crack growth with a slight delay.

\paragraph{Crack Propagation Behavior.}
Under the same boundary and domain conditions, we then examined the evolution of the damage variable. As shown in Figures~\ref{2DED2-1} and \ref{2DED2-2}, the crack in the standard DF-PFM propagated horizontally, driven by the compressive loading—even with the presence of the P-wave. However, in the modified model with the unilateral contact condition, the crack propagated more naturally along the initial angled direction. The timing of the crack propagation in the DF-PFM with a unilateral contact condition is relatively later than the one in the DF-PFM.
These distinctions arise because the standard DF-PFM includes a compression term in the elastic energy density, which can lead to unphysical overlapping deformation and horizontal crack propagation. In contrast, the modified model excludes this term, resulting in more realistic shear-driven crack behavior under compression.

\section{Conclusion}

In this study, we proposed and analyzed two dynamic fracture phase field models (F-PFMs) for simulating fault-type earthquake ruptures. The first model \eqref{DF}, referred to as the dynamic F-PFM (DF-PFM), extends the classical quasi-static fracture formulation to incorporate elastic wave propagation. The second model \eqref{DFU}, the DF-PFM with a unilateral contact condition, introduces one-sided contact constraints to prevent nonphysical interpenetration along the crack surfaces. For both models, we established energy dissipation identities under appropriate regularity assumptions, ensuring the energy-consistency of the formulations.

To support numerical implementation, we developed linear implicit time-discrete schemes and derived their corresponding weak formulations. These schemes were implemented using P1 finite element methods, and some numerical experiments were conducted to validate and compare the models.

A key simulation involved a subsurface fault configuration: an obliquely oriented initial crack under compressive loading, impacted by a downward-traveling P-wave from above. Both models reproduced fault rupture triggered by the incident wave. However, the resulting crack propagation patterns exhibited qualitative differences. The DF-PFM, lacking contact constraints, produced nonphysical kink-type crack paths associated with negative opening displacements. In contrast, the DF-PFM with a unilateral contact condition generated more realistic shear-driven rupture behavior consistent with fault-slip mechanisms.

These findings demonstrate that incorporating contact constraints can significantly improve the physical fidelity of dynamic fracture simulations in geophysical contexts. We believe that the proposed DF-PFM with a unilateral contact condition can also be extended to model other dynamic fracture processes, such as fluid-driven fracture, thermal cracking, and desiccation-induced damage in geomaterials.

Looking ahead, several promising directions remain for future research. From a modeling and computational standpoint, incorporating more realistic physical ingredients—such as rate-and-state friction laws, nonlinear material responses, or thermo-poroelastic couplings—could enhance the predictive power of the framework. In addition, the use of adaptive mesh refinement and high-order time integration schemes may improve computational efficiency and resolution, particularly for large-scale three-dimensional simulations.

Beyond fault dynamics, the proposed DF-PFM with a unilateral contact condition can be extended to a wide range of fracture phenomena in geomechanics and material science. Coupling this framework with multiphysics processes would enable the study of complex fracture interactions in both natural and engineered systems.


\section*{Acknowledgements}
Md Mamun Miah was supported by the Japanese Government Monbukagakusho (MEXT) Scholarship in the years 2022-2025. Masato Kimura was partially supported by JSPS KAKENHI Grant Numbers JP24H00184, JP25K00920.

\end{document}